\documentclass[11pt]{elsarticle}

\usepackage{charter} % font
\usepackage{etex}
\usepackage{ifpdf}
\usepackage{amsmath} 
\usepackage{enumerate}
\usepackage{booktabs,siunitx}
\usepackage{graphicx,subfigure,epstopdf,epsfig}
\usepackage[hang,flushmargin]{footmisc} 
\usepackage{sansmath}
\usepackage{color,soul}
\sethlcolor{yellow}
\usepackage{xcolor}
\usepackage{hyperref}
\usepackage{amssymb}                  
\usepackage[all]{xy}
\usepackage[width=16cm,height=24cm,top=2cm]{geometry}
\usepackage{algorithm,algpseudocode} %package install: sudo apt-get install texlive-science
\algnewcommand\algorithmicinput{\textbf{Input:}}
\algnewcommand\Input{\item[\algorithmicinput]}
\algnewcommand\algorithmicoutput{\textbf{Output:}}
\algnewcommand\Output{\item[\algorithmicoutput]}
\usepackage{multirow}
\usepackage[font=small,labelfont=bf]{caption}
\newcommand{\norm}[1]{\left\lVert#1\right\rVert}

\DeclareMathOperator*{\argmin}{argmin}

\def\KKK{{K}}
\def\SSS{{\mathcal S}}
\def\VVV{{\bf V}}

\usepackage{lineno}
\usepackage{fullpage}

\journal{Journal of Computational Physics}

\begin{document}

\begin{frontmatter}

 \title{Fast Spherical Centroidal Voronoi Mesh Generation: A Lloyd-preconditioned LBFGS Method in Parallel}
 \tnotetext[t1]{H. Yang and M. Gunzburger partially supported by US Department of Energy Office of Science grant number DE-SC0016591. L. Ju  partially supported by US Department of Energy Office of Science grant number DE-SC0016540. }
    
  \author[stu]{Huanhuan Yang}
  \ead{huan2yang@stu.edu.cn}

  \author[fsu]{Max Gunzburger}
  \ead{mgunzburger@fsu.edu}

  \author[usc]{Lili Ju}
  %\corref{mycorrespondingauthor}}
  \ead{ju@math.sc.edu}

  \address[stu]{Department of Mathematics, Shantou University, Shantou, Guangdong 515063, China}
  \address[fsu]{Department of Scientific Computing, Florida State University, Tallahassee, FL 32306, USA}
  \address[usc]{Department of Mathematics, University of South Carolina, Columbia, SC 29208, USA}

\begin{abstract}
Centroidal Voronoi tessellation (CVT)-based mesh generation is a very effective technique for creating high-quality Voronoi meshes
 and their dual Delaunay triangulations that  often {play} a crucial role in applications, including ocean and atmospheric  simulations using finite volume schemes. In the next generation climate 
 models, the spacing scales change dramatically across the whole sphere and require ultra-high resolution and smooth transitions from coarse to fine grid regions. Thus fast and robust
 spherical CVT (SCVT) meshing algorithms become highly desirable.
In this paper, we first propose a Lloyd-preconditioned limited-memory BFGS method for constructing SCVTs that is also applicable to 
{the construction} of CVTs  of general domains. This method  is then  parallelized {based on} overlapping domain decomposition, enabling excellent scalability on 
distributed systems.  Results of several computational experiments show that the new method could incur computational time costs one order of magnitude smaller  
compared  with some  existing methods for generating large-scale highly variable-resolution meshes, while also providing significantly improvements in mesh quality.
\end{abstract}


\end{frontmatter}

%%%%%%%%%%%%%%%%%%%%%%%%%%%%%%%%%%%%%%%%%%%%%%%%%%%%%%%%%%%%%%%%%%%%%%%%%%%
%%%%%%%%%%%%%%%%%%%%%%%%%%%%%%%%%%%%%%%%%%%%%%%%%%%%%%%%%%%%%%%%%%%%%%%%%%%
%%%%%%%%%%%%%%%%%%%%%%%%%%%%%%%%%%%%%%%%%%%%%%%%%%%%%%%%%%%%%%%%%%%%%%%%%%%
\section{Introduction}\label{intro}

A centroidal Voronoi tessellation (CVT) \cite{du1999centroidal} of a given domain is a Voronoi tessellation whose generators coincide with the mass centroids of their corresponding Voronoi cells. Because of its ability to construct point sets that are locally quasi-uniformly distributed with respect to a given point-density function, CVT is a desirable technique for creating high-quality Voronoi and Delaunay meshes \cite{du2002grid, du2003constrained}. For example, such meshes are often crucial for global or regional climate modeling using finite volume schemes\cite{ringler2011exploring,ringler2013multi,JAME20337,JAME20336}. 

CVTs can be constructed using probabilistic methods such as MacQueen's method \cite{macqueen1967some} or probabilistic Lloyd methods \cite{ju2002probabilistic} that are elegant random sequential sampling methods that do not require the explicit construction of Voronoi tessellations. However, these probabilistic methods feature very slow convergence so that more efficient deterministic methods for CVT construction are more often used. Typical deterministic CVT construction algorithms include Lloyd's method \cite{lloyd1982least, du1999centroidal}, the Lloyd-Newton method \cite{DuAccel2006}, and quasi-Newton methods such as limited-memory BFGS (LBFGS) method \cite{Liu2009CVT,Hateley2015}. 

In addition to the geometric definition given above, CVTs also have an analytical definition as being a critical point of an ``energy'' functional associated with a set of generators 
\cite{du1999centroidal}. In particular, a local minimizer of a CVT energy functional determines a stable CVT. Lloyd's method can be viewed as a linearly convergent gradient descent optimization method \cite{du1999centroidal} with respect to the CVT energy function. The Lloyd-Newton method is an attempt to accelerate CVT construction by taking advantage of its quadratic rate of convergence. However, it is not often used for high-resolution CVT grid generation due to the burdensome computation of a Hessian inverse at each iteration. 
 
Quasi-Newton methods are usually more feasible choices because they make use of efficient approximate Hessian matrices or their inverses, and generally have super-linear convergence. In \cite{Liu2009CVT} the LBFGS quasi-Newton method and a preconditioned version thereof are used to construct CVTs in planar regions or surfaces. The preconditioner used in \cite{Liu2009CVT} is the incomplete Cholesky factorization of a modified Hessian matrix; however, this approach is not guaranteed to be stable or effective \cite{benzi1999comparative, Hateley2015} and, in addition, it has been shown that such LBFGS methods preconditioned in this manner do not apparently improve the efficiency of LBFGS \cite{Liu2009CVT} in terms of computational time.  A graph Laplacian preconditioner is used in \cite{Hateley2015} which reduces the computational time to 56.9\% on average for constructing CVT meshes with certain non-uniform point-density functions. However, the graph Laplacian preconditioner can be troublesome for generating high-resolution meshes not only because of the limited efficiency due to having to repeatedly solve a large-scale linear system, but also because it is limited effectiveness in reducing the CVT energy.

In recent years, the need for improving the efficiency of the next generation of the computational climate models has driven the need for the more efficient construction of large-scale meshes on the sphere, meshes that feature high-quality and highly variable-resolution with smooth transitions. In this paper, we propose a fast and effective quasi-Newton method for the construction of spherical CVTs (SCVT) but which is also applicable to the construction of CVTs in more general domains. One of our key contributions is the development of a Lloyd-preconditioned LBFGS algorithm; specifically, the Lloyd step is iteratively executed and used as approximations of varying initial Hessian inverses. For quasi-uniform meshes, our Lloyd-preconditioned LBFGS scheme performs similarly to LBFGS. However, for highly variable-resolution meshes,  it can dramatically speed up LBFGS while also providing significantly improved mesh quality. Also differing from all previous efforts using LBFGS for computing CVTs, we develop a parallel implementation of our method on distributed systems by using an overlapping domain decomposition approach \cite{JacobMPI2013}. Our parallel algorithm  features well-balanced loading of mesh points and has excellent performance with respect to  strong scaling efficiency. Parallel test results show almost linear speedup when the number of mesh points is of the order of millions.

The remainder of the paper is organized as follows. In Section \ref{prelim}, we review the definitions and basic properties of CVTs and their extension to surfaces/manifolds and also review some popular CVT construction algorithms. Then, in Section \ref{cvt-alg}, we propose the Lloyd-preconditioned LBFGS method for computing SCVTs and its parallelization on distributed systems. Several numerical experiments and comparisons are presented in Section \ref{results} to demonstrate  the performance of our method in both serial and parallel environments. Some concluding remarks are made in Section \ref{conclusion}.

%%%%%%%%%%%%%%%%%%%%%%%%%%%%%%%%%%%%%%%%%%%%%%%%%%%%%%%%%%%%%%%%%%%%%%%%%%%
%%%%%%%%%%%%%%%%%%%%%%%%%%%%%%%%%%%%%%%%%%%%%%%%%%%%%%%%%%%%%%%%%%%%%%%%%%%
%%%%%%%%%%%%%%%%%%%%%%%%%%%%%%%%%%%%%%%%%%%%%%%%%%%%%%%%%%%%%%%%%%%%%%%%%%%
\section{Review of centroidal Voronoi tessellation and its generalizations}\label{prelim}

In this section, we first introduce some concepts and properties of centroidal Voronoi tessellation and then review some existing algorithms for their construction.

%%%%%%%%%%%%%%%%%%%%%%%%
%%%%%%%%%%%%%%%%%%%%%%%%
\subsection{Definitions and properties of CVTs}
%%%%%%%%%%%%%%%%%%%%%%%%
%%%%%%%%%%%%%%%%%%%%%%%%

Given a bounded domain $\Omega \subset \mathbb{R}^N$ and a set of points $\{\mathbf{z}_i\}_{i=1}^\KKK \subset \overline{\Omega}$, the  {\it Voronoi region} $V_i$ corresponding to the {\it generator} (or {\it site}) $\mathbf{z}_i$ is defined as
\begin{equation}\label{viInR}
V_i = \{\mathbf{x} \in \Omega \,\,\,:\,\,\, \|\mathbf{x} - \mathbf{z}_i\| < \|\mathbf{x} - \mathbf{z}_j\| \,\,\,\mbox{ for}\,\,\, j=1,\ldots,\KKK \,\, \,\mbox{and}\,\,\, j \neq i\},
\end{equation}
where $\|\cdot\|$ denotes the Euclidean norm in $\mathbb{R}^N$. A {\it Voronoi tessellation} or {\em Voronoi diagram} $\VVV = \{V_i\}_{i=1}^\KKK$ forms a special partition of $\Omega$ in the sense that $V_i \cap V_j=\emptyset$ for $i\neq j$ and $\overline{\Omega} = \cup_i\overline{V_i}$. Voronoi cells are convex polytopes except possibly for those whose boundary intersect with the boundary of $\Omega$. The dual graph of a Voronoi tessellation corresponds to the Delaunay triangulation of the set of Voronoi generators.

Given a point-density function $\rho(\mathbf{x})>0$ defined over $\overline{\Omega}$, the mass center $\mathbf{z}^*_i$ of a Voronoi cell $V_i$ is defined by
\begin{equation}\label{masscentroid}
\mathbf{z}^*_i = \dfrac{\displaystyle\int_{V_i}\mathbf{y}\rho(\mathbf{y})d\mathbf{y}}{\displaystyle\int_{V_i}\rho(\mathbf{y})d\mathbf{y}}, \qquad i = 1, \ldots, \KKK. 
\end{equation}
In general, $\mathbf{z}_i \neq \mathbf{z}^*_i$.  In the special case that $\mathbf{z}_i = \mathbf{z}^*_i$ for $i=1,\ldots,\KKK$, i.e., each generator coincides with the mass centroid of the corresponding Voronoi cell,  the Voronoi tessellation is referred to as a {\it centroidal Voronoi tessellation} (CVT) \cite{du1999centroidal}.

Alternately and equivalently to the just given geometric characterization of CVTs, one can define them through a variational approach. Letting $\mathbf{Z} = \{\mathbf{z}_i\}_{i=1}^\KKK$ denote an arbitrary ordered sequence of points in $\Omega$ and $\VVV=\{V_i\}_{i=1}^\KKK$ an arbitrary tessellation of $\Omega$, we define the ``energy'' functional
$$
\mathcal{F}(\mathbf{Z},\VVV) =  \sum\limits_{i=1}^\KKK\int_{V_i}|\mathbf{y}-\mathbf{z}_i|^2\rho(\mathbf{y})d\mathbf{y}.
$$
A CVT is a critical point of this functional \cite{du1999centroidal} in the sense that among all possible sets of $\KKK$ points $\mathbf{Z}$ in $\Omega$ and among all possible tessellations $\VVV$ of $\Omega$ into $\KKK$ nonoverlapping subregions, this functional is rendered stationary if and only if $\{\mathbf{Z},\VVV\}$ form a CVT. Furthermore, local minimizers $\{\mathbf{Z},\VVV\}$ of $\mathcal{F}$ determine stable CVTs. Stable CVTs are more desirable in practical applications than unstable ones which correspond to saddle points of $\mathcal{F}$.

Next, consider the tessellation of a compact surface or manifold $\SSS \subset \mathbb{R}^N$. Given a set of points  $\{\mathbf{z}_i\}_{i=1}^\KKK \subset \SSS$, we can define the corresponding Voronoi tessellation of $\SSS$ similarly to (\ref{viInR}) as
\begin{equation}\label{viInS}
V^c_i = \{\mathbf{x} \in \SSS \,\,\,:\,\,\, \|\mathbf{x} - \mathbf{z}_i\| < \|\mathbf{x} - \mathbf{z}_j\| \,\,\mbox{ for }\,\, j=1,\ldots,\KKK, \,\, j \neq i\},
\end{equation}
For each Voronoi cell $V^c_i$, we use a generalized definition for the center of mass to ensure that it lies on the surface $\SSS$. Specifically, we refer to $\mathbf{z}^c_i$ as the {\it constrained mass centroid} \cite{du2003constrained} of $V^c_i$ on $\SSS$ if 
\begin{equation}\label{ccentroid}
\mathbf{z}^c_i = \argmin_{\mathbf{z}\in\SSS}\int_{V^c_i}\|\mathbf{y}-\mathbf{z}\|^2\rho(\mathbf{y})d\sigma(\mathbf{y}),
\end{equation}
where the point-density function $\rho(\mathbf{x})>0$ is defined on the surface $\SSS$ and $d\sigma(\cdot)$ denotes the surface area element. Note that in general the classical mass centroid $\mathbf{z}^*_i$ defined by \eqref{masscentroid} does not lie on $\SSS$. However, as pointed out in \cite{du2003constrained}, there is a connection between the constrained mass centroid $\mathbf{z}^c_i$ and the classical mass centroid $\mathbf{z}^*_i$, namely, if $\mathbf{z}^c_i \in \SSS-\partial\SSS$, then $\mathbf{z}^*_i - \mathbf{z}^c_i$ is normal to the surface $\SSS$ at $\mathbf{z}^c_i$, i.e., $\mathbf{z}^c_i$ is the projection of $\mathbf{z}^*_i$ onto $\SSS$ along the normal direction to $\SSS$ at $\mathbf{z}^c_i$. This feature is especially useful if $\SSS$ is a sphere because on the one hand, the construction of a constrained mass centroid using its definition \eqref{ccentroid} is considerably more difficult compared to the construction of the classical mass centroid and, on the other hand, for $\SSS$ a sphere, the former is merely the intersection of the radial line going through the latter with the sphere.

The tessellation of $\SSS$ defined by (\ref{viInS}) is referred as a {\it constrained centroidal Voronoi tessellation} (CCVT) if $\mathbf{z}_i = \mathbf{z}^c_i$ for $i=1,\ldots,\KKK$, i.e., the generator associated with each Voronoi cell $V^c_i$ coincides with the constrained mass centroid of that cell. It is also proved in \cite{du2003constrained} that a stable CCVT can also be defined as a local minimizer of the CCVT energy function
$$
\mathcal{F}^c(\mathbf{Z},\VVV)  = \sum\limits_{i=1}^\KKK\int_{V^c_i}\|\mathbf{y}-\mathbf{z}_i\|^2\rho(\mathbf{y})d
\sigma(\mathbf{y}).
$$

It is worth emphasizing that although the generators of CCVTs are constrained to lie on a surface/manifold, the distance metric used is still the Euclidean distance rather than the geodesic distance. 

When $\SSS$ is a sphere in $\mathbb{R}^3$ or a subset of a sphere, as is the case in ocean modeling, a CCVT is referred to as {\em spherical} CVT (SCVT).

%%%%%%%%%%%%%%%%%%%%%%%%
%%%%%%%%%%%%%%%%%%%%%%%%
\subsection{Some existing deterministic CVT construction algorithms}
%%%%%%%%%%%%%%%%%%%%%%%%
%%%%%%%%%%%%%%%%%%%%%%%%

\subsubsection{Lloyd's method}

An elegant deterministic method for computing CVTs/CCVTs is Lloyd's method which is a natural byproduct of the geometric characterization of CVTs. Starting from an initial placement of $\KKK$ points in $\Omega$, a Lloyd iteration determines the Voronoi diagram of $\Omega$ associated with those points, then moves each point to the mass centroid of its Voronoi cell. This pair of steps is repeated until a sufficient small decrease in the point movement is met. Lloyd's method can be viewed as a fixed-point iteration for the Lloyd map
$$
 \mathbf{T}: ~\mathbf{Z} \to \left\{ \frac{\displaystyle\int_{V_i(\mathbf{Z})}\mathbf{y}\rho(\mathbf{y})d\mathbf{y}}{\displaystyle\int_{V_i(\mathbf{Z})}\rho(\mathbf{y})d\mathbf{y}} \right\}_{i=1}^\KKK,
$$
where $\{V_i(\mathbf{Z})\}_{i=1}^\KKK$ denotes the Voronoi tessellation associated with points $\{\mathbf{z}_i\}_{i=1}^\KKK$. Convergence analyses for some typical cases are given in \cite{du1999centroidal, du2006convergence}. 

Based on the variational characterization of a CVT as a critical point of the CVT energy function $\mathcal{F}(\mathbf{Z},\VVV)$, a CVT can also be constructed by derivative-based optimization methods because the energy function $\mathcal{F}$ is $C^2$ smooth for a convex domain \cite{Liu2009CVT}. Recognizing that in Lloyd's method $\VVV=\VVV(\mathbf{Z})$, i.e., the tessellation $\VVV$ is the Voronoi tessellation corresponding to $\mathbf{Z}$, we have that $\mathcal{F}=\mathcal{F}(\mathbf{Z},\VVV(\mathbf{Z}))$, i.e., we can view $\mathcal{F}$ as a function of $\mathbf{Z}$ only. Then the gradient of $\mathcal{F}(\mathbf{Z},\VVV(\mathbf{Z}))$ is given by
$$
\frac{d\mathcal{F}}{d\mathbf{z}_i} = \frac{\partial \mathcal{F}}{\partial\mathbf{z}_i}(\mathbf{Z}, \VVV(\mathbf{Z})) = 2\mathbf{z}_i\int_{V_i}\rho(\mathbf{y})d\mathbf{y} - 2\int_{V_i}\mathbf{y}\rho(\mathbf{y})d\mathbf{y},
$$
where we use the fact $\frac{\partial \mathcal{F}}{\partial \VVV} = 0$ at $(\mathbf{Z}, \VVV(\mathbf{Z}))$ due to the property of Voronoi regions \cite{iri1984fast, du1999centroidal}. Denoting the mass $\displaystyle\int_{V_i}\rho(\mathbf{y})d\mathbf{y}$ by $m_i$, we have
$$
m_i^{-1}\frac{d\mathcal{F}}{d\mathbf{z}_i}  = 2\mathbf{z}_i - 2\mathbf{T}_i(\mathbf{Z}).
$$ 
Thus, the Lloyd iteration $\mathbf{Z}^{(n+1)} = \mathbf{T}(\mathbf{Z}^{(n)})$ can be rewritten as
\begin{equation}\label{loy-itr}
{\mathbf{Z}^{(n+1)} = \mathbf{Z}^{(n)} - \big(2\mathbf{M}^{(n)}\big)^{-1} \nabla\mathcal{F}(\mathbf{Z}^{(n)})},
\end{equation}
with the mass matrix $\mathbf{M}^{(n)}$ being a diagonal matrix $\mbox{diag}(m_i)$ associated with the Voronoi tessellation $\VVV(\mathbf{Z}^{(n)})$.
In the context of optimization, Lloyd's method is a special case of a gradient descent method. Although it has only a linear rate of convergence, the Lloyd iteration always decreases the energy function $\mathcal{F}$ until a local minimizer is obtained without the need of line searches for step-size control.

\subsubsection{Newton and  quasi-Newton methods}

An attempt to accelerate CVT construction is the Lloyd-Newton method \cite{DuAccel2006}. Although  the Newton iteration has quadratic rate of convergence, it is usually not feasible for large scale CVT grid generation due to the burden of storing and computation of Hessian inverses. For this reason, quasi-Newton methods are better choices because they use approximations of the Hessian matrices or their inverses and have super-linear rates of convergence. Specifically, given an approximation $\mathbf{B}^{(n)}$ of the Hessian matrix $\mathbf{H}(\mathbf{Z}^{(n)})$ of $\mathcal{F}(\mathbf{Z}^{(n)})$, a quasi-Newton iteration operates as follows: for $n=0,1,2,\ldots$
\begin{enumerate}
\item[1:] solve $\mathbf{B}^{(n)}\mathbf{q} = - \nabla\mathcal{F}(\mathbf{Z}^{(n)})$;
\item[2:] update $\mathbf{Z}^{(n+1)} = \mathbf{Z}^{(n)} + \alpha^{(n)} \mathbf{q}$.
\end{enumerate}
Here $\mathbf{q}$ is a search direction along which an optimal step size $\alpha^{(n)}$ is determined by a line search algorithm to satisfy the Wolfe conditions \cite{Nocedal2006NO}. The choice $\mathbf{B}^{(n)} = \mathbf{H}(\mathbf{Z}^{(n)})$ results in the classical Newton method. 

One of the most popular quasi-Newton methods is the BFGS method \cite{Nocedal2006NO}. In the standard BFGS approach, at the $n$-th step, the inverse of Hessian matrix is iteratively updated by a symmetric matrix $\widetilde{\mathbf{H}}^{(n+1)}$ being close to the last Hessian inverse approximation and satisfying the secant equation 
$$\widetilde{\mathbf{H}}^{(n+1)} y_n = s_n,\quad \mbox{where } s_{n}=\mathbf{Z}^{(n+1)}-\mathbf{Z}^{(n)}, ~ y_n = \nabla\mathcal{F}(\mathbf{Z}^{(n+1)}) - \nabla\mathcal{F}(\mathbf{Z}^{(n)}).$$
Such derivation forces the quadratic forms of function $\mathcal{F}$ around $\mathbf{Z}^{(n)}$ and $\mathbf{Z}^{(n+1)}$ match gradients respectively at these two steps. 
The BFGS updating formula in the product form can be given explicitly as
\begin{eqnarray*}
\widetilde{\mathbf{H}}^{(n)} & = &\big(S^T_{n-1}\cdots S^T_{n-m(n)}\big)\widetilde{\mathbf{H}}^{(0)}\big(S_{n-m(n)}\cdots S_{n-1}\big) \\
&& +\big(S^T_{n-1}\cdots S^T_{n-m(n)+1}\big)\rho_{n-m(n)}s_{n-m(n)}s_{n-m(n)}^T\big(S_{n-m(n)+1}\cdots S_{n-1}\big) \\
&& +\big(S^T_{n-1}\cdots S^T_{n-m(n)+2}\big)\rho_{n-m(n)+1}s_{n-m(n)+1}s_{n-m(n)+1}^T\big(S_{n-m(n)+2}\cdots S_{n-1}\big) \\
&&+ \cdots + \rho_{n-1}s_{n-1}s_{n-1}^T,
\end{eqnarray*}
where 
$$\rho_i = 1/(y_i^Ts_i), \qquad S_i = \mathbf{I}-\rho_{i}s_{i}s_{i}^T.$$ 
In the classic BFGS updating scheme, $m(n)=n$ is taken and the initial approximation $\widetilde{\mathbf{H}}^{(0)}$ is often set to be a constant multiple of the identity matrix.

In the limited-memory BFGS (LBFGS), one {needs to} specify $M$, the number  of BFGS corrections needed to be stored. That is, only the recent $M$ pairs of $\{s_i, y_i\}$ are kept. The BFGS approximation is modified by setting $m(n) = \min\{n, M\}$ and allowing the initial approximation $\widetilde{\mathbf{H}}^{(0)}$ to vary from iteration to iteration. For instance, $\widetilde{\mathbf{H}}^{(0)}$ is typically replaced by $\widetilde{\mathbf{H}}^{(0)}_n = \gamma_n\mathbf{I}$, where
\begin{equation*}
\gamma_n = \frac{y_{n-1}^Ts_{n-1}}{y_{n-1}^Ty_{n-1}}
\end{equation*}
is close to the reciprocal of an eigenvalue of $\mathbf{H}(\mathbf{Z}^{(n)})$.
In practical computation, the product $\widetilde{\mathbf{H}}^{(n)}\nabla\mathcal{F}(\mathbf{Z}^{(n)})$, rather than the
matrix $\widetilde{\mathbf{H}}^{(n)}$, is recursively updated. Specifically,  the LBFGS updating at step $n$ is done as follows:

\vspace{0.15cm}

\begin{algorithmic}[1]\em
%\Input 
%\Output 
\State Initialize $ \mathbf{q} \leftarrow -\nabla\mathcal{F}(\mathbf{Z}^{(n)})$;
%\Statex \% backward loop:
\For{$i=n-1, \cdots, n-m(n)$} \hspace{0.2cm} \% backward loop
	\Statex \hspace{0.5cm}$\alpha_i \leftarrow \rho_is_i^T\mathbf{q}$;
	\Statex \hspace{0.5cm}$\mathbf{q} \leftarrow \mathbf{q} - \alpha_iy_i$;
\EndFor

\vspace{0.15cm}

\State $\boxed{\mathbf{q} \leftarrow \widetilde{\mathbf{H}}^{(0)}_n\mathbf{q};}$\label{varInitH}

\vspace{0.15cm}

\For{$i=n-m(n), \cdots, n-1$} \hspace{0.2cm} \% forward loop
	\Statex \hspace{0.5cm}$\mathbf{q} \leftarrow \mathbf{q} + s_i(\alpha_i-\rho_iy_i^T\mathbf{q})$;
\EndFor
\State Update $\mathbf{Z}^{(n+1)} = \mathbf{Z}^{(n)} + \alpha^{(n)} \mathbf{q}$. \hspace{0.02cm} \% line search
\end{algorithmic}

\vspace{0.15cm}

%%%%%%%%%%%%%%%%%%%%%%%%%%%%%%%%%%%%%%%%%%%%%%%%%%%%%%%%%%%%%%%%%%%%%%%%%%%
%%%%%%%%%%%%%%%%%%%%%%%%%%%%%%%%%%%%%%%%%%%%%%%%%%%%%%%%%%%%%%%%%%%%%%%%%%%
%%%%%%%%%%%%%%%%%%%%%%%%%%%%%%%%%%%%%%%%%%%%%%%%%%%%%%%%%%%%%%%%%%%%%%%%%%%
\section{Lloyd-preconditioned LBFGS in parallel}\label{cvt-alg}

In this section, we will first develop a Lloyd preconditioner for the LBFGS method for SCVT construction  so that its performance can be further improved, after which we propose an efficient parallelization algorithm for the method on distributed systems based on overlapping domain decomposition, 
which is important in terms of both computational and storage efficiency. 

%%%%%%%%%%%%%%%%%%%%%%%%
%%%%%%%%%%%%%%%%%%%%%%%%
\subsection{Lloyd-preconditioned LBFGS}
%%%%%%%%%%%%%%%%%%%%%%%%
%%%%%%%%%%%%%%%%%%%%%%%%

The idea of a preconditioned LBFGS method (P-LBFGS) proposed in \cite{jiang2004preconditioned} is to choose a matrix $\mathbf{B}^{(n)}$ very alike to the Hessian $\mathbf{H}(\mathbf{Z}^{(n)})$ and periodically replace line \ref{varInitH} in the LBFGS iteration by
\begin{equation}\label{precStep}
\mathbf{q} \leftarrow (\mathbf{B}^{(n)})^{-1} \mathbf{q}.
\end{equation}
It requires that $\mathbf{B}^{(n)}$ should be inexpensively achievable, and most importantly, the above equation can be easily solved.  If  the matrix $\mathbf{B}^{(n)}$ can capture significant features of the Hessian, then the convergence of LBFGS can be hastened.

It has been shown that the CVT energy function has second-order smoothness for convex domains \cite{Liu2009CVT} and that the Hessian matrix can be computed exactly; see \cite{iri1984fast,Liu2009CVT}. Denote by $J_i$ the set of node indices of the Voronoi cells adjacent to $V_i$ (the Voronoi cell corresponding to generator $\mathbf{z}_i$).  Let the coordinates of the point $\mathbf{z}_i$ (resp.~$\mathbf{y}$) be denoted by $z_{ik}$ (resp.~$y_k$), $k=1, \ldots, N$. Then, the second-order derivatives of the CVT energy are given by the following explicit formulas:  for $k,l=1, \ldots, N$,
\begin{equation}\label{hessForm}
\left\{ 
\begin{array}{ll}
\frac{\partial^2 \mathcal{F}}{\partial z_{ik}^2} = 
2 m_i - \displaystyle\sum_{j\in J_i}\displaystyle \int_{\overline{V_i}\cap \overline{V_j}}\frac{2}{\norm{\mathbf{z}_j - \mathbf{z}_i}}(z_{ik}-y_\KKK)^2\rho(\mathbf{y})d{\bf y}, &
\\[0.3cm]
\frac{\partial^2 \mathcal{F}}{\partial z_{ik}\partial z_{il}} = 
-\displaystyle\sum_{j\in J_i} \displaystyle\int_{\overline{V_i}\cap \overline{V_j}}\frac{2}{\norm{\mathbf{z}_j - \mathbf{z}_i}}(z_{ik}-y_\KKK)(z_{il}-y_l)\rho(\mathbf{y})d{\bf y}, & 
k\neq l, \\[0.3cm]
\frac{\partial^2 \mathcal{F}}{\partial z_{ik}\partial z_{jl}} = 
\displaystyle\int_{\overline{V_i}\cap \overline{V_j}}\frac{2}{\norm{\mathbf{z}_j - \mathbf{z}_i}}(z_{ik}-y_\KKK)(z_{jl}-y_l)\rho(\mathbf{y})d{\bf y}, & j\in J_i, \\[0.3cm]
\frac{\partial^2 \mathcal{F}}{\partial z_{ik}\partial z_{jl}} = 0, & j\neq i, ~j\notin J_i.
\end{array}
\right.
\end{equation}
In \cite{Liu2009CVT}, a modified Cholesky factorization of the Hessian matrix $\mathbf{H}(\mathbf{Z}^{(n)})$ is applied to obtain the positive-definite matrix $\mathbf{B}^{(n)}$, and perform (\ref{precStep}) every $\widetilde N$ iterative steps, where $\widetilde N$ is a user-defined parameter. The numerical experiments in \cite{Liu2009CVT} show that LBFGS and P-LBFGS yield better performance compared to Lloyd's method and Newton's method. However, the preconditioned version of LBFGS used in \cite{Liu2009CVT} only improves LBFGS a little, and these two perform very similarly even in the case of non-uniform point-density functions.

Another approach at developing a P-LBFGS method is to use a graph Laplacian preconditioner (graph Laplacian P-LBFGS), proposed in \cite{Hateley2015}. They consider CVT construction in a {polygonal} domain $\Omega$ in two dimensions. Let $T_{ij}$ denote the triangle inside $V_i$ containing the point $\mathbf{z}_i$ and the edge $
\overline{V_i}\cap \overline{V_j}$. If $V_i$ has an edge $e_{ib}$ in $\partial \Omega$, denote the triangle formed by $e_{ib}$ and $\mathbf{z}_i$ as $T_{ib}$. Then, the graph Laplacian matrix $\mathbf{A} = (a_{ij})$ is given by
$$
\left\{ 
\begin{array}{ll}
a_{ii} =\displaystyle \sum_{j\in J_i}\int_{T_{ij}\cup T_{ji}}\rho(\mathbf{y})d\mathbf{y} + 2\sum_{e_{ib}} \int_{T_{ib}}\rho(\mathbf{y})d\mathbf{y}, \\[0.3cm]
a_{ij} = \displaystyle-\int_{T_{ij}\cup T_{ji}}\rho(\mathbf{y})d\mathbf{y}, &\quad j\in J_i,\\[0.3cm]
a_{ij} = 0,	& \quad j\neq i, ~j\notin J_i.
\end{array}
\right.
$$
According to the numerical tests in \cite{Hateley2015}, for a constant point-density function, the graph Laplacian P-LBFGS has similar efficiency compared to the classical LBFGS algorithm, but for non-uniform point-density functions, the  graph Laplacian preconditioner can reduce the computational time by 56.9\% on average. 
However, the use of graph Laplacian preconditioners to speed up LBFGS could be limited for ocean grids due to some difficulties: first, solving the Laplacian linear system on large-scale {meshes} at each iteration is computationally expensive;  second, the Laplacian matrix $\mathbf{A}$ is singular for closed surfaces such as the surface of the sphere and thus special care is needed in choosing an appropriate solution for the preconditioning process. Furthermore, the graph Laplacian preconditioner sometimes {does not work} well in significantly reducing the CVT energy, as we reported in Section ~\ref{GLvsLP}.

In this paper, we propose a new preconditioner (initial Hessian inverse, line \ref{varInitH} in the LBFGS iteration), inspired by Lloyd's iteration, as
\begin{equation}\label{loyInit}
\mathbf{\widetilde{H}}^{(0)}_n = \Big(2\mathbf{M}^{(n)}\Big)^{-1} = \mbox{diag}\bigg(2\int_{V_i(\mathbf{Z}^{(n)})}\rho(\mathbf{y})d\mathbf{y}\bigg)^{-1}.
\end{equation}
We call this the {\it Lloyd-preconditioned LBFGS} (Lloyd P-LBFGS) method because this approach is equivalent to taking the initial BFGS direction to be the Lloyd's updating direction (see (\ref{loy-itr})) and the initial line-search step size as 1. 

There are several reasons behind this idea. First of all, because the Lloyd map $\mathbf{T}$ is continuous, the Lloyd iteration converges globally to a critical point of $\mathcal{F}$ if the iterations stay in a compact set \cite{du2006convergence}. Although the compactness of the iteration has not been rigorously justified in the 
literature, it seems to be intuitively true and matches practical experiences. This nice property of Lloyd's method is not shared by Newton or quasi-Newton methods. Second, the initial Hessian inverse given by (\ref{loyInit}) is very simple while still capturing a lot of Hessian matrix information. Note that $2\mathbf{M}^{(n)}$ is the main diagonal part of the Hessian matrix as shown in (\ref{hessForm}). A connection of (\ref{loyInit}) to the graph Laplacian matrix is more obvious: it is an {approximation} of 
the diagonal components of the Laplacian matrix $\mathbf{A}$ because 
$$\int_{V_i} = \sum_{j\in J_i}\int_{T_{ij}} + \displaystyle\int_{T_{ib}}
$$ 
and we can treat $
\displaystyle
\int_{T_{ij}\cup T_{ji}}$ approximately as $2\displaystyle\int_{T_{ij}}$. Lastly, the diagonal matrix $\mathbf{M}^{(n)}$ itself contains significant information about the density distribution. 

To apply LBFGS for minimizing the energy function $\mathcal{F}^c$ in the case of CCVT, in \cite{Liu2009CVT} a ``pseudo'' derivative is taken by projecting the gradient of the CVT energy onto the tangent plane of the compact surface $\SSS$. That is, they proposed to use the tangential derivative
$$
\left.\frac{\partial \mathcal{F}^c}{\partial \mathbf{z}_i}\right|_{\SSS} = \frac{\partial \mathcal{F}^c}{\partial \mathbf{z}_i} - \bigg(\frac{\partial \mathcal{F}^c}{\partial 
\mathbf{z}_i}\cdot \mathbf{N}(\mathbf{z}_i)\bigg) \mathbf{N}(\mathbf{z}_i)
$$
in the minimization process, where $\mathbf{N}(\mathbf{z}_i)$ is the normal to the surface $\SSS$ at $\mathbf{z}_i$. The updated set of generators needs to be projected back onto the surface at each iteration.
Following this approach, we can derive the Lloyd preconditioner for {the case} of the unit sphere  analogous to that for planar domains. Note that 
$$
\frac{\partial \mathcal{F}^c}{\partial \mathbf{z}_i} = 2\mathbf{z}_i\int_{V_i^c}\rho(\mathbf{y})d\sigma - 2\int_{V_i^c}\mathbf{y}\rho(\mathbf{y})d\sigma.
$$ 
Setting $\mathbf{c}_i = \displaystyle\int_{V_i^c}\mathbf{y}\rho(\mathbf{y})d\sigma$, the tangential derivative can be written as
$$
\left.\frac{\partial \mathcal{F}^c}{\partial \mathbf{z}_i}\right|_{\SSS} = 2(\mathbf{c}_i^T\mathbf{z}_i)\mathbf{z}_i - 2\mathbf{c}_i
$$
so that we have
$$
\mathbf{c}_i = (\mathbf{c}_i^T\mathbf{z}_i)\mathbf{d}_i, \quad\mbox{ with } \mathbf{d}_i = \mathbf{z}_i - (2\mathbf{c}_i^T\mathbf{z}_i)^{-1} \left.\frac{\partial \mathcal{F}
^c}{\partial \mathbf{z}_i}\right|_{\SSS}.
$$
The updating by $\mathbf{d}_i$ is then a Lloyd iteration because the projection of $\mathbf{d}_i$ onto the unit surface coincides with that of $\mathbf{c}_i$. As such, we can take the Lloyd preconditioner in $\mathbb{R}^{K\times K}$ as
$$
\mathbf{\widetilde{H}}^{(0)}_n = \mbox{diag}\bigg(2(\mathbf{c}_i^{(n)})^T\mathbf{z}_i^{(n)}\bigg)^{-1}.
$$

%An alternate way for generating SCVT is to perform the optimization process in a two-dimensional
%space rather than three-dimensional, by using polar coordinates. Consider the unit sphere in $\mathbb{R}^3$. According to the spherical coordinate mapping
%$$
%\mathbf{p} = [p_a, p_o]^T \longrightarrow \mathbf{z} 
%= [\cos(p_a)\cos(p_o), \cos(p_a)\sin(p_o), \sin(p_a)]^T,
%$$
%the exact gradient of the energy functional $\mathcal{F}^c$ can be computed by
%$$
%\frac{\partial \mathcal{F}^c}{\partial \mathbf{p}_i} = \frac{\partial \mathcal{F}^c}{\partial \mathbf{z}_i}\frac{\mathcal{D} \mathbf{z}_i}{\mathcal{D} \mathbf{p}_i},
%$$
%where $\{\mathbf{p}_i\}_{i=1}^\KKK$ are the generators of the spherical Voronoi regions expressed in polar coordinates. The Lloyd iteration is now explicitly given by (using Taylor's expansion)
%$$
%\mathbf{z}^{(n+1)}_i - \mathbf{z}^{(n)}_i = \mathbf{J}_i^{(n)} \big(\mathbf{p}^{(n+1)}_i - \mathbf{p}^{(n)}_i\big),
%$$
%where $\mathbf{J}_i^{(n)} = \left.\frac{\mathcal{D}\mathbf{z}}{\mathcal{D}\mathbf{p}}\right|_{\mathbf{p}_i^{(n)}}$ is the Jacobian,
%so that we can approximate the Lloyd step for updating $\mathbf{p}_i$ as
%$$
%\mathbf{p}^{(n+1)}_i = \mathbf{p}^{(n)}_i- \Big(2m_i^{(n)}(\mathbf{J}_i^{(n)})^T\mathbf{J}_i^{(n)}\Big)^{-1}
%\left.\frac{\partial\mathcal{F}^c}{\partial\mathbf{p}_i}\right|_{\mathbf{p}_i^{(n)}}.
%$$
% The Lloyd preconditioner  in the $\mathbb{R}^{2K\times 2K}$ form is finally given by
%$$
%\widetilde{\mathbf{H}}^{(0)}_n = \mbox{diag}\bigg(2m_i^{(n)}(\mathbf{J}_i^{(n)})^T\mathbf{J}_i^{(n)} \bigg)^{-1}.
%$$

%%%%%%%%%%%%%%%%%%%%%%%
%%%%%%%%%%%%%%%%%%%%%%%
\subsection{Parallel implementation}
%%%%%%%%%%%%%%%%%%%%%%%
%%%%%%%%%%%%%%%%%%%%%%%
Our parallel implementation of the Lloyd-preconditioned LBFGS method for SCVT mesh generation is based on overlapping domain decomposition.

%%%%%%%%%%%%%%%%%%%%%%%
\subsubsection{Parallelization of Delaunay triangulations on sphere}
%%%%%%%%%%%%%%%%%%%%%%%

As is true for all other deterministic methods for CVT/SCVT construction, the Lloyd-preconditioned LBFGS algorithm involves explicit integrations over Voronoi cells. The Voronoi cells can be determined from the dual Delaunay triangulation. In this subsection we present a parallel algorithm for computing Delaunay triangulations and spherical Delaunay triangulations first studied in \cite{JacobThesis, JacobMPI2013}.

\paragraph{Planar Delaunay triangulations} 
A Delaunay triangulation for a set of points $\{\mathbf{z}_i\}_{i=1}^\KKK$ in a plane is a triangulation $\mathcal{T}$ such that no point in $\{\mathbf{z}_i\}_{i=1}^\KKK$ is strictly inside the circumcircle of any triangle in $\mathcal{T}$. Delaunay triangulations maximize the minimum angle of all the angles of the triangles in the triangulation; they tend to avoid extremely acute angles of the triangles. The uniqueness of Delaunay triangulation for $\{\mathbf{z}_i\}_{i=1}^\KKK$ is guaranteed if the set of points is in general position (no set of four points lie on a circle whose interior does not contain points from $\{\mathbf{z}_i\}_{i=1}^\KKK$).

\begin{figure}[!h]
\begin{center}
\includegraphics[scale=0.5]{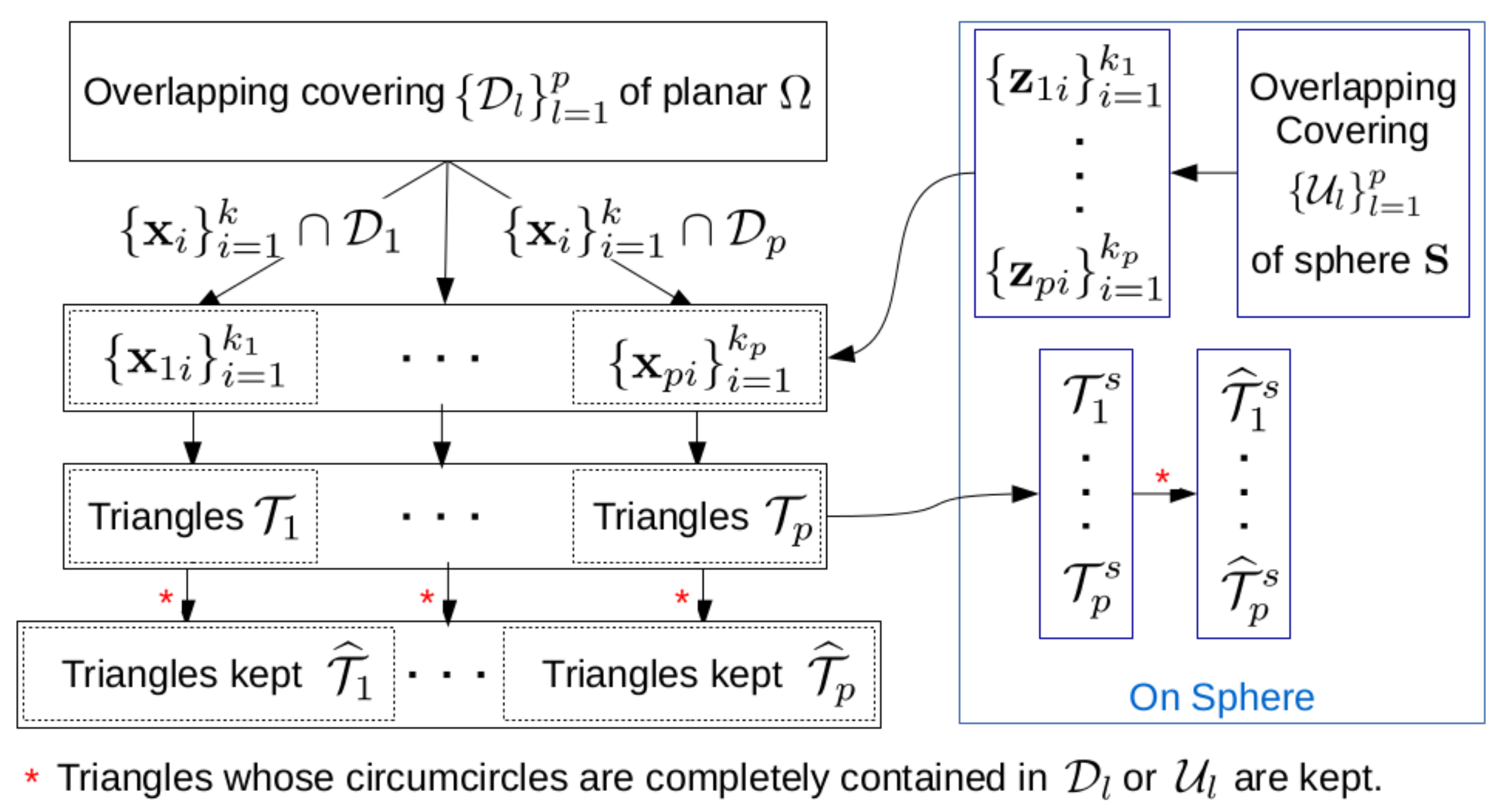}
\caption{Parallelization of the Delaunay triangulation. The left part is for a planar domain, the right is for the sphere. Stereographic projections are used from right to left.}
\label{paraDT}
\end{center}
\end{figure}

The parallelization of planar Delaunay triangulation algorithms has been widely studied in the past decades; see, e.g., \cite{cignoni1998dewall, dymond2001}. Popular algorithms commonly split the global point set into equally distributed subsets, each of which can then be triangulated simultaneously in parallel. The resulting local triangulations are then stitched together to form a global triangulation. Such a merging step may need to modify significant portions of the local triangulations, thus it is the most significant step that affects time efficiency. In this work, we resort to an alternative {merging} step that needs no modifications of the local triangulations and thus performs well in parallel. The algorithm for planar Delaunay triangulation consists of the following three steps (see Figure~\ref{paraDT}-left).
\begin{enumerate}
\item Given a planar domain $\Omega$ and a set {of} points $\{\mathbf{z}_i \}_{i=1}^\KKK$ in $\Omega$, we preform an overlapping covering of $\Omega$ and record the overlapping connectivity or neighboring list. A simple example of the covering $\mathcal{C}$ can be a set of $p$ disks denoted by $\{\mathcal{D}_l(\mathbf{p}_l, r_l)\}_{l=1}^p$, where $\mathbf{p}_l$ and $r_l$ are the disk center and radius respectively. The $\KKK$ points in $\Omega$ are then distributed into $p$ subregions, each of which is handled by a single processor. 
Because of the overlapping feature of the covering, some points may belong to multiple processors.
\item On the $l$-th processor, a local Delaunay triangulation $\mathcal{T}_l$ is carried out independently. To this end, one can employ a serial Delaunay algorithm such as the Delaunay triangulator in the Triangle software package \cite{Shewchuk1996}. 
\item 
On each processor, we only keep triangles whose circumcircles are  completely contained inside the corresponding subregion. Denoting the new triangulation on the $l$-th processor by $\widehat{\mathcal{T}}_l$, then the circumcircle of any triangle in $\cup_l\widehat{\mathcal{T}}_l$ does not 
contain any points from $\{\mathbf{z}_i \}_{i=1}^\KKK$.
{Each local triangulation $\widehat{\mathcal{T}}_l$ is then exactly a portion of the global Delaunay triangulation due to  the uniqueness of Delaunay triangulations. }
\end{enumerate}

\paragraph{Spherical Delaunay triangulation} 
We now consider spherical Delaunay triangulation on the unit sphere $\SSS$ in $\mathbb{R}^3$.
In this case, triangles and Voronoi cells are replaced by spherical ones whose boundaries are geodesic arcs. The planar parallelization scheme is combined with stereographic projections.

A stereographic projection is a conformal mapping that projects a sphere onto a plane tangent to {the sphere. 
Letting $\mathbf{t}$ denote} the contact point of the sphere and its tangent plane, the stereographic projection of a point $\mathbf{z}$ on the unit sphere to the point $
\mathbf{x}$ on the tangent plane is given by (in Cartesian coordinates)
$$ 
\mathcal{P}: \mathbf{z} \to \mathbf{x} = s\mathbf{z} + (s-1)\mathbf{t}, \quad \mbox{ where } s = \dfrac{2}{(\mathbf{t}\cdot (\mathbf{z}+\mathbf{t}))}.
$$
The conformality of stereographic projections implies the preservation of circumcircles of triangles along with their interiors and hence the Delaunay criteria. Therefore, the Delaunay triangulation of a subregion of the sphere can be obtained by the stereographic mapping from the planar Delaunay tessellation on the corresponding tangent plane. The adaptation to spherical Delaunay triangulation from a planar one is described by the following three steps; see Figure~\ref{paraDT}.
\begin{enumerate}
\item We preform an overlapping covering of the unit sphere $\SSS$ and record the overlapping connectivity or neighboring list. The covering $\mathcal{U}$ could simply be a set of geodesic disks $\{\mathcal{U}_l(\mathbf{t}_l, r_l)\}_{l=1}^p$, where $\mathbf{t}_l$ is the disk center and $r_l$ is the geodesic radius. Points $\{\mathbf{z}_i\}_{i=1}^\KKK$ on $\SSS$ are then distributed into $p$ subregions with overlaps.
\item On the $l$-th processor, we project the point set $\{\mathbf{z}_i\}_{i=1}^\KKK \cap \mathcal{U}_l(\mathbf{t}_l, r_l)$ onto the plane tangent to the sphere at $\mathbf{t}_l$. We then proceed to construct a planar Delaunay triangulation $\mathcal{T}^p_l$. The local spherical Delaunay triangulation $\mathcal{T}^s_l$ is then produced by stereographic projection from $\mathcal{T}^p_l$. 
\item In each subregion, say $\mathcal{U}_l(\mathbf{t}_l, r_l)$, we only keep spherical triangles whose circumcircles are completely contained inside, so as to guarantee the Delaunay property with respect to the global point set $\{\mathbf{z}_i\}_{i=1}^\KKK$. In particular, for disk subregions, let $\mathbf{p}^c_j$ and $r^c_j$ denote the center and radius of the circumcircle of the $j$-th spherical triangle in $\mathcal{U}_l(\mathbf{t}_l, r_l)$; the triangle selection criteria is $\arccos(\mathbf{t}_l \cdot \mathbf{p}^c_j) + r^c_j \leq r_l.$ Each local spherical triangulation is now exactly a portion of the global spherical Delaunay triangulation.
\end{enumerate}

To obtain optimal load balancing and robustness, the target domain should be decomposed based on the given point-density function. A precomputed ultra-coarse CVT is then a natural choice for effecting the partition. The main challenge is now the determination of the sizes of the overlaps. The overlapping regions should be large enough to make sure no true Delaunay triangles are missing in the final triangulation. A simple overlapping decomposition for quasi-uniform meshes is to use disk regions centered at partition cell centers with radius being the maximum distance from the center to its adjacent centers. More challenges exist for variable resolution meshes. In this work, we take the Voronoi sort method, i.e., the union of the owned partition cell and its immediately adjacent partition cells defines the overlapping domain decomposition for parallel triangulation. Failure may occur when the global point set does not include enough points for the decomposition of choice. A practical suggestion on the number of grid points is at least 16 multiples of the number of partitions, or a number corresponding to at least two-level bisections from the partition CVT, as pointed out in \cite{JacobMPI2013}.

%%%%%%%%%%%%%%%%%%%%%%%
\subsubsection{Parallelization of Lloyd P-LBFGS}
%%%%%%%%%%%%%%%%%%%%%%%

\begin{figure}[!h]
\begin{center}
\includegraphics[scale=0.5]{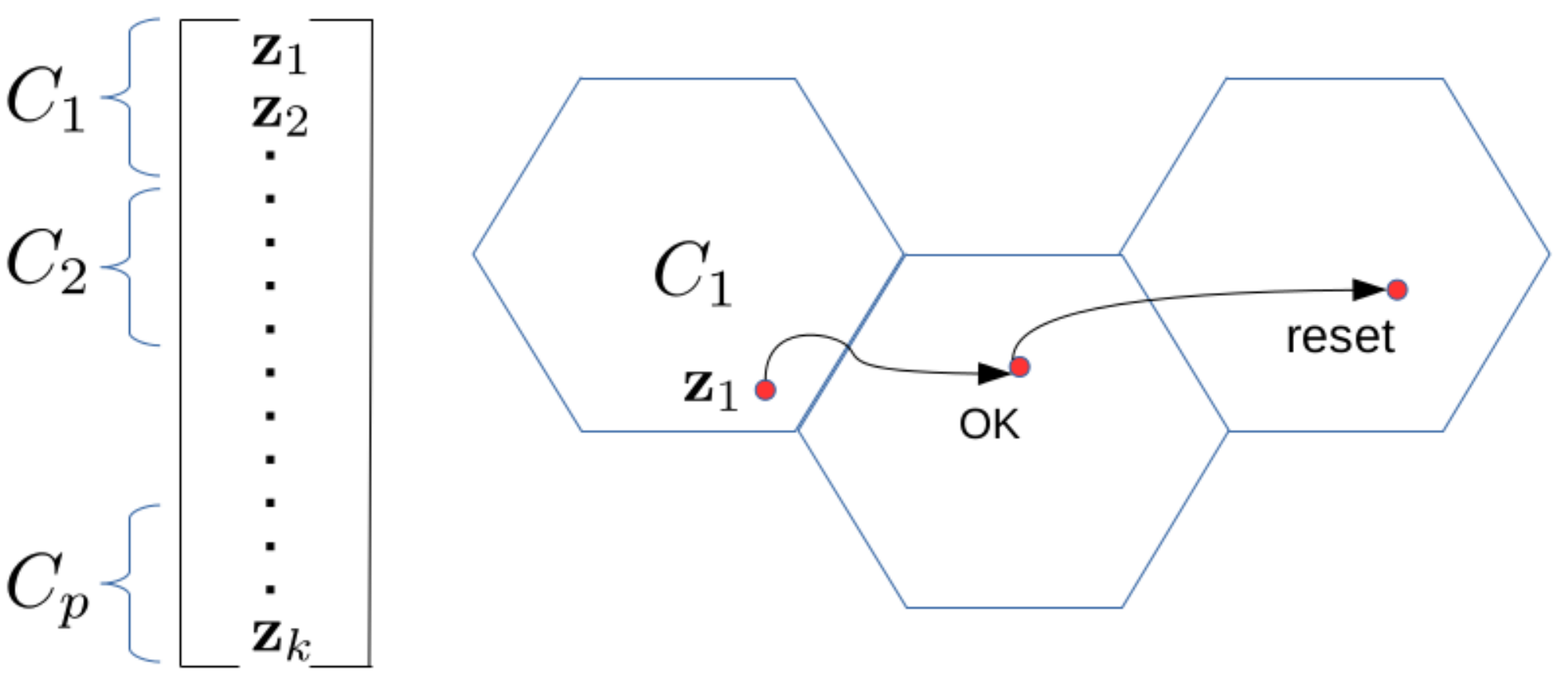}
\caption{The global point vector is initially decomposed by non-overlapping partitions. We assume Voronoi points move at most into their neighboring partition cells; if a point moves further, reset LBFGS.}
\label{partBFGS}
\end{center}
\end{figure}

To parallelize the Lloyd P-LBFGS scheme, we first transform the global point vector $\mathbf{Z} = (\mathbf{z}_i)_{i=1}^\KKK$ into non-overlapping local vectors, so that arithmetic operations like additions and dot products can be efficiently applied on the distributed vectors. This non-overlapping decomposition can be initially carried out according to {the disjoint partition cells $C_1, \ldots, C_p$; see} Figure~\ref{partBFGS}-left. On each processor, the global IDs of points in the distributed vector are recorded and will be used for vector assembly after point updating. For any point, say $\mathbf{z}_i$, if it is initially in partition cell  $l$ by  disjoint partition, we call cell $l$ its {\it arithmetic} position. As the iterations proceed, $\mathbf{z}_i$ could move from one partition cell  to another, and we call the current partition cell it belongs to as its {\it geometric} position. Whereas the geometric position could vary, we always fix the arithmetic position for convenience of computations inside LBFGS.

The evaluations of the energy function, the gradient, and the Lloyd preconditioner require integrations over Voronoi cells. All integrations are actually computed based on the dual Delaunay triangulation. Given a triangle $\triangle(\mathbf{z}_a, \mathbf{z}_b, \mathbf{z}_c)$ in the Delaunay tessellation, with its circumcenter $\mathbf{c}$ and three edge midpoints $\mathbf{z}_{ab}, \mathbf{z}_{bc}$, and $\mathbf{z}_{ca}$, the triangle can be split into six sub-triangles, {each made} of an original triangle vertex, an edge midpoint, and the circumcenter. The sub-triangle $\triangle(\mathbf{z}_a, \mathbf{z}_{ab}, \mathbf{c})$, for instance, is part of the Voronoi cell whose generator is located at $\mathbf{z}_a$. Should the orientation of sub-triangles be taken into account,  an integration over a Voronoi cell is the sum of integration  over such sub-triangles inside the cell. 

The {evaluation} of the functional gradient and the Lloyd preconditioner described above are performed in 
parallel and stored according to the geometric positions of grid points. To be clear, we will call the partition cell associated with which current processor is handling as ``my cell'', and the union with its overlapped layers as ``my region''. After points are assigned into ``my region'', local Delaunay triangulation is performed. We then integrate over ``my cell'' for energy function evaluation. In the meanwhile, at each point in ``my cell'', we compute the function gradient and the Lloyd preconditioner. Non-blocking communication is then used to assemble their arithmetic vectors, whose distribution is the same as the non-overlapping partition of the initial global points.

It is worth noting that in our code there is no need for global gathering of local Delaunay triangulations or updated mesh points (Voronoi generators). The communications for updating mesh points are only among two-level neighboring partition cells, under the assumption that grid points move at most into their neighboring partition cells. Here the ``two-level'' approach is to ensure the updating of points in overlapped layers. The assumption is generally fulfilled in practice by consecutively performing optimization-bisection until the final level of optimization; this is detailed in next section. As such, the communications for assembling the gradient and preconditioner are only in one level, i.e.,~only among neighboring partition cells. When the assumption is not satisfied, a reset on P-LBFGS is executed (see the illustration in Figure~\ref{partBFGS}-right), i.e., the arithmetic position for each mesh point is re-computed.

%%%%%%%%%%%%%%%%%%%%%%%
%%%%%%%%%%%%%%%%%%%%%%%
\subsection{Description of the parallel algorithm}\label{algSec}
%%%%%%%%%%%%%%%%%%%%%%%
%%%%%%%%%%%%%%%%%%%%%%%

There are still several ingredients that need to be considered for parallel P-LBFGS for SCVT construction, such as what are the point-density function, the initialization, and the stopping criteria. 

The point-density function $\rho(\mathbf{x})$ plays a crucial role in how the converged mesh points are distributed. For the numerical tests in this paper, we use analytic point-density functions that feature significant physically relevant information drawn from climate modeling. A point-density function of 
course also can  be formulated from a given initial mesh, or even from a local mesh of a subregion of the sphere. 

The convergence of iterative methods is often highly sensitive to the initial configuration. Exploring a good initialization for SCVT computation has been an interesting topic in the field. To the best of our knowledge, a Monte Carlo initialization based on the point-density function is often suitable for coarse mesh generation. For high-resolution meshes, however, a consecutive optimization-bisection procedure is almost mandatory to obtain an initial guess that is good enough to save computational cost. 

The stopping criteria for P-LBFGS here are based on the maximum number of iterations, the norm of gradient of the energy function, and the step size of mesh point movement. 

In summary, the parallel version of the Lloyd  P-LBFGS algorithm for SCVT  construction consists of the following steps.
	
\vspace{0.15cm}

\begin{algorithmic}[1]\em
%\Input 
%\Output 
\State Build the partition information from a coarse CVT;
\State Set $n \leftarrow 0$;
\State Decompose $\mathbf{Z}^{(n)}$ into disjoint partition cells accordingly; store global IDs of  points; \label{initialPart}
\While{stopping criteria not satisfied} 
	\State Preprocess:\label{preproc}
	\Statex	\hspace{1cm} assign points from $\mathbf{Z}^{(n)}$ into ``my region'' (2-level neighboring comm.~);
	\Statex \hspace{1cm} go to Step \ref{initialPart} if a point in $\mathbf{Z}^{(n)}$ moves out of its initial neighboring partition cells
	\Statex	\hspace{1cm} for points in ``my region'', build local Delaunay triangulation;
	\State For points in ``my cell'', compute $\mathcal{F}(\mathbf{Z}^{(n)})$, $\nabla\mathcal{F}(\mathbf{Z}^{(n)})$, $\widetilde{\mathbf{H}}_n^{(0)}$;
	\State Sum $\mathcal{F}(\mathbf{Z}^{(n)})$; assemble $\nabla\mathcal{F}(\mathbf{Z}^{(n)})$, $\widetilde{\mathbf{H}}_n^{(0)}$ to arithmetic positions (1-level 
	neighboring comm.~)
	\State Compute search direction $\mathbf{q}$ by P-LBFGS;
	\State Update $\mathbf{Z}^{(n+1)} = \mathbf{Z}^{(n)} + \alpha^{(n)} \mathbf{q}$ by line search: 
	\Statex \hspace{1cm} occasionally perform Step \ref{preproc} on $\mathbf{Z}^{(n)} + s \mathbf{q}$ and evaluate $\mathcal{F}$, for some step size $s$.
	\State $n \leftarrow n+1$
\EndWhile
\State Gather all updated points;
\State Build local Delaunay triangulation and then gather all triangles (non-repeatedly).
\end{algorithmic}

%%%%%%%%%%%%%%%%%%%%%%%%%%%%%%%%%%%%%%%%%%%%%%%%%%%%%%%%%%%%%%%%%%%%%%%%%%%
%%%%%%%%%%%%%%%%%%%%%%%%%%%%%%%%%%%%%%%%%%%%%%%%%%%%%%%%%%%%%%%%%%%%%%%%%%%
%%%%%%%%%%%%%%%%%%%%%%%%%%%%%%%%%%%%%%%%%%%%%%%%%%%%%%%%%%%%%%%%%%%%%%%%%%%
\section{Performance tests}\label{results}

%In this section, we will first compare the Lloyd-preconditioned LBFGS method with the Lloyd's method and the LBFGS method to demonstrate the %improvement on 
%computational efficiency and mesh quality. Then, we compare the Lloyd preconditioner with the graph Laplacian preconditioner. After that, %the strong scaling 
%efficiency of the parallel P-LBFGS solver is studied and we will show numerical results of comparing three iterations in parallel.

Our code for the Lloyd-preconditioned LBFGS solver for SCVT construction is written in C++ and MPI; the Delaunay triangulation is handled by the ``Triangle'' package, the LBFGS routine  is from the serial HLBFGS package by Yang Liu (in which the line-search routine is converted from the original LBFGS Fortran code by Jorge Nocedal). We set the maximum number of function evaluations for the line search to ten and store at most seven BFGS corrections in memory ($M=7$), as suggested in \cite{Liu2009CVT}. The SCVT grid generator in this study is under active development and the code is available at https://github.com/hyang52/scvt\_HLBFGS.

The quality of the Voronoi cell $V$ of the generated meshes is quantified by  $cellQ$ 
$$
 cellQ(V) = \frac{\min \mbox{cell edge}}{\max \mbox{cell edge}}
$$
and the triangle quality $triQ$, which is defined on an triangle $\triangle$ with edge lengths $a, b, c$ by
\begin{equation}\label{Qdef}
 triQ(\triangle) = \frac{(a+b-c)(b+c-a)(c+a-b)}{abc}.
 \end{equation}

\begin{figure}[!h]
\begin{center}
\includegraphics[width=2in]{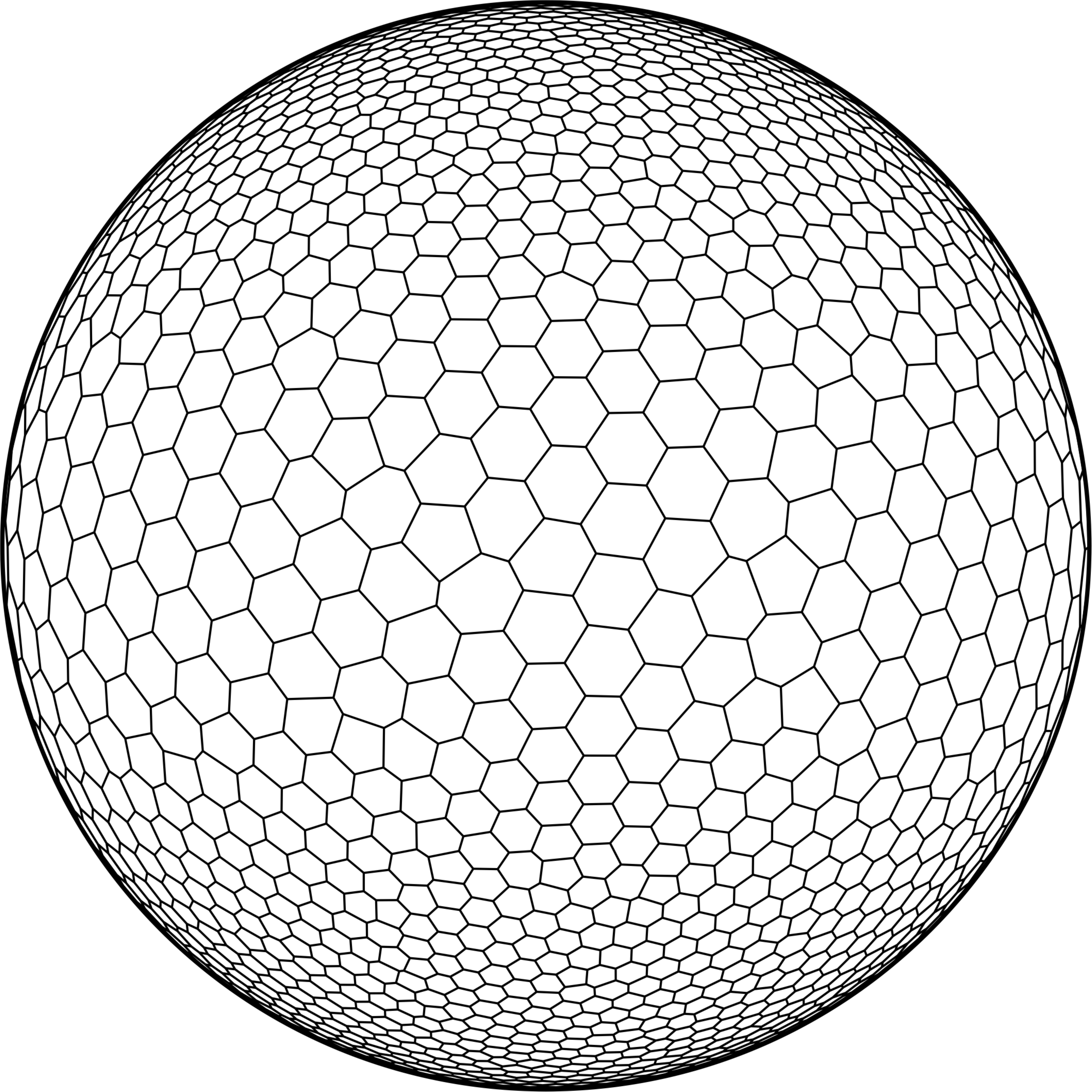}
\hspace{0.0cm}
\includegraphics[width=2in]{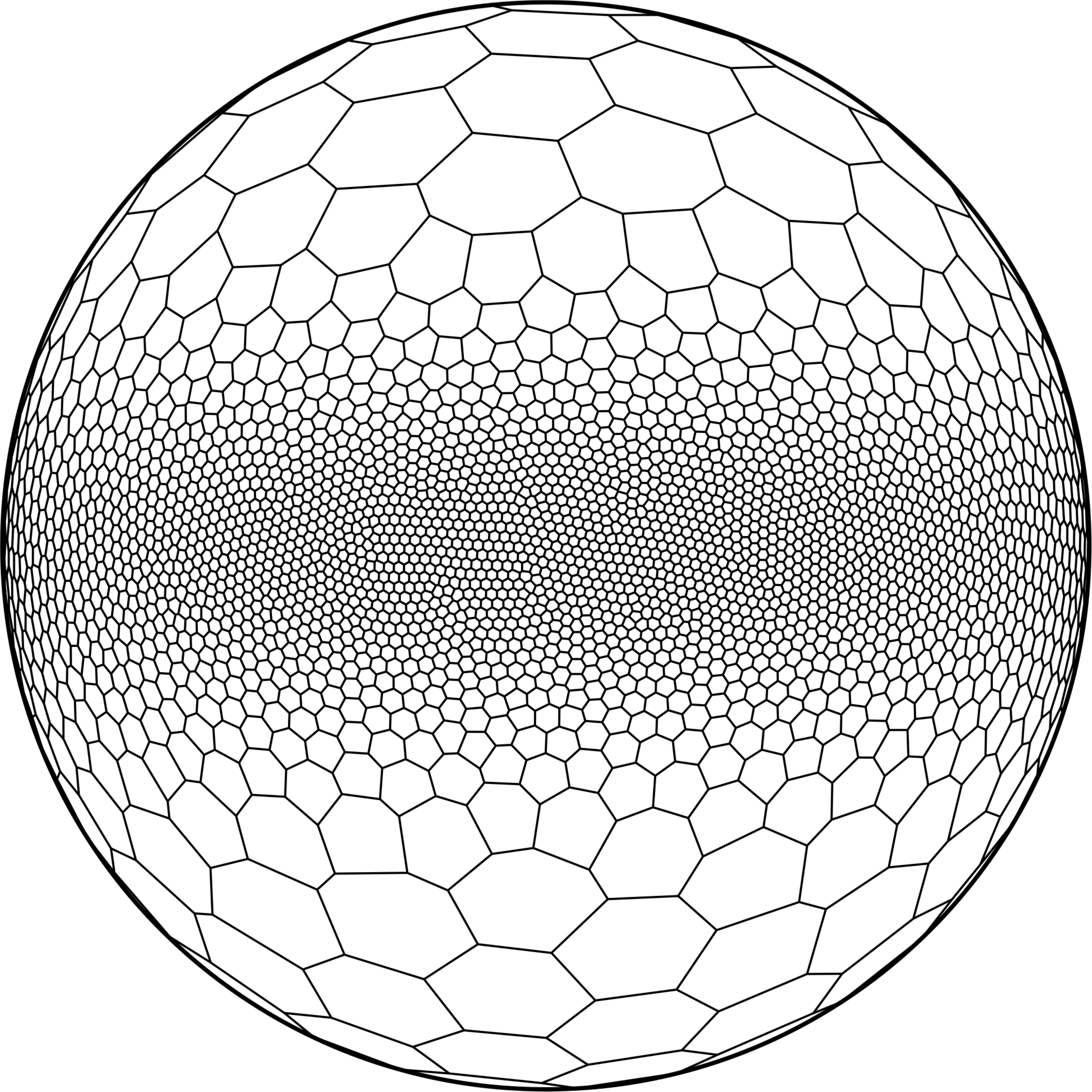}
\hspace{0.0cm}
\includegraphics[width=2in]{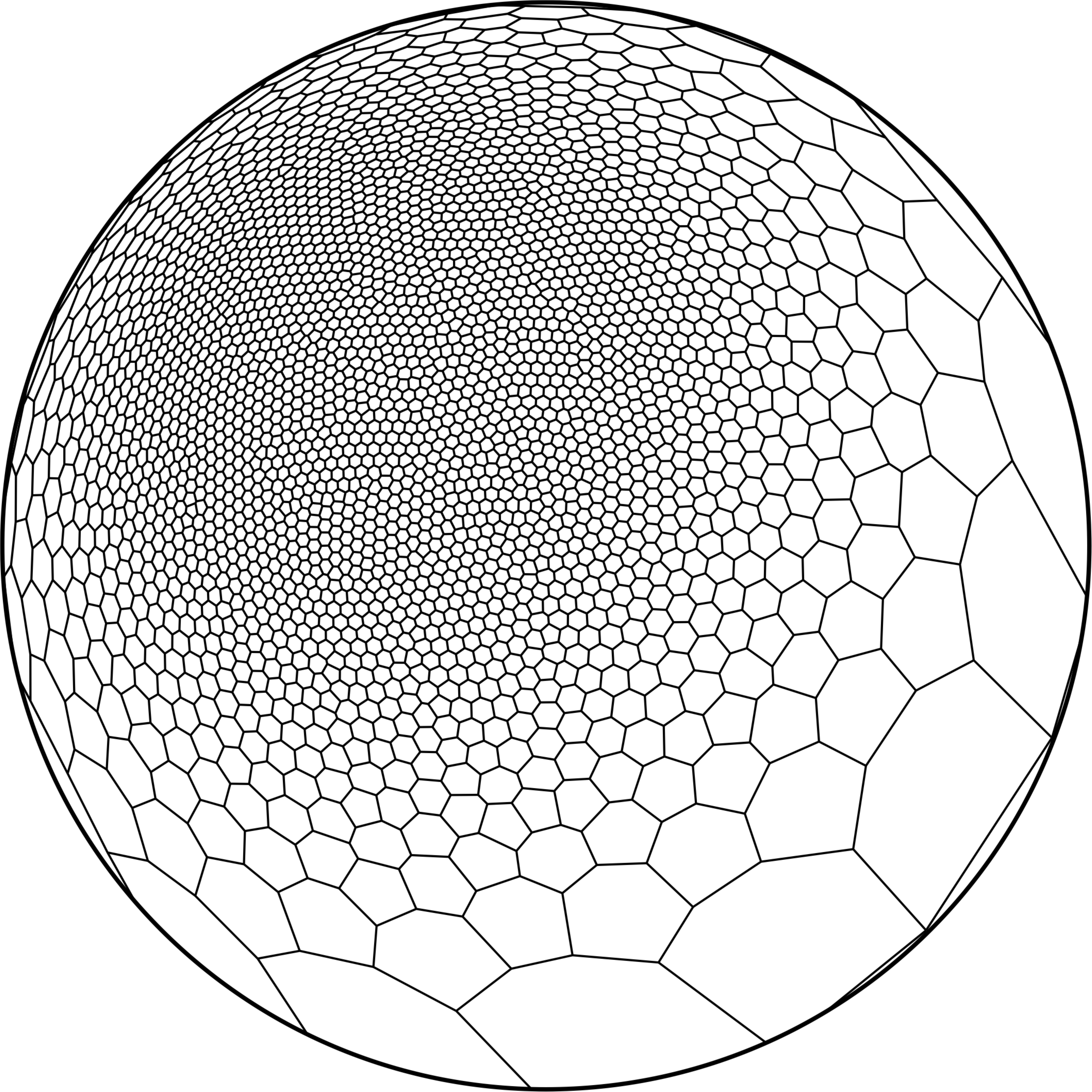}
\end{center}
\caption{The SCVT meshes with 2,562 generators created by Lloyd P-LBFGS with three point-density functions: X3, X16, X64 (from left to right).}
\label{VTcomp}
\end{figure}

Three different point-density functions defined on the unit sphere are used in our tests to represent some typical cases. For each point-density function, we use the max-min ratio of the point-density to control the max-min ratio of mesh size because it is conjectured that for any two Voronoi cells $V_i$ and $V_j$ from a SCVT, their mesh spacings are related to the point-density function as
$$ \frac{h_i}{h_j} \approx \bigg(\frac{\rho(\mathbf{z}_j)}{\rho(\mathbf{z}_i)} \bigg)^{1/4}.$$
The first point-density function we used is given by 
$$
\rho(x,y,z) = (1-\gamma)z^4 + \gamma
$$ 
with $\gamma = (1/3)^4$, so that the mesh size on the equator is three times of that on the poles. This size distribution is typically used in some simulations of ocean modeling. We will call it the X3 point-density; see the SCVT mesh in the left column of Figure~\ref{VTcomp}. 

The other two point-density functions are defined as
$$
\rho(\mathbf{z}_i) = \frac{1}{2(1-\gamma)}\bigg[\tanh\Big(\frac{\beta-d(\mathbf{z}_i)}{\alpha}\Big) +1 \bigg] + \gamma,
$$
where $d(\mathbf{z}_i)$ is a distance function that can be simply the geodesic distance from $\mathbf{z}_i$ to a pre-specified fixed point $\mathbf{z}_c$ on the sphere. This definition results in a relatively large value of $\rho$ within a distance $\beta$ from the point $\mathbf{z}_c$. The point-density transitions from a large value to small value across a geodesic distance of $\alpha$. To have a variety of tested densities, we define $d(\mathbf{z}_i)$ in different ways for the second and third point-density functions. Specifically, for the second point-density function, we define 
$$
d(\mathbf{z}_i) = \bigg(\frac{|\mathbf{z}_i - \mathbf{z}_{ic}|_{\SSS}^2}{w_{\rm lon}^2}  + \frac{|\mathbf{z}_i - \mathbf{z}_{ci}|_{\SSS}^2}{w_{\rm lat}^2}\bigg)^{1/2},
$$
where $\mathbf{z}_{ic}$ is a point whose latitude and longitude coincide with $\mathbf{z}_i$ and  $\mathbf{z}_c$, respectively, and vice versa for $\mathbf{z}_{ci}$. The norm $|\cdot|_{\SSS}$ denotes geodesic distance. We choose the widths $w_{\rm lon}$ and $w_{\rm lat}$ to be $0.3$ and $1.2$ with {$\mathbf{z}_c=(1.0,0.0,0.0)$}, and set $\beta = \pi/6$ and $\alpha = 0.3$. The minimum point-density is set to be $\gamma = (1/16)^4$, so we refer to the defined function as the X16 point-density; see the SCVT grid in the middle column of Figure~\ref{VTcomp}). 

For the third point-density, we take a function used in \cite{ringler2011exploring} for a multi-resolution modeling with the shallow water equations. The distance is now $d(\mathbf{z}_i) = |\mathbf{z}_i - \mathbf{z}_c|_{\SSS}$ with {$\mathbf{z}_c=(0.0, -0.866, 0.5)$}, the parameters are taken as $\beta = \pi/6$ and $\alpha = 0.15$. Here, we set the minimum point-density $\gamma = (1/64)^4$ and refer to the function as the X64 point-density; see the SCVT grid in the right column  of Figure~\ref{VTcomp}).
There are two reasons of using this X64 point-density. On the one hand, it corresponds to a highly variable multi-resolution mesh that is needed in global ocean and coastal system modeling for which the mesh size may range from a $O(10^4)$km global scale to $O(10^{-1})$km local scale. On the other hand, it can be used in global spherical mesh generation for the purpose for limited-area climate modeling, if one culls the coarse portion from the global one.
 
%%\begin{figure}[!h]
%%\begin{center}
%%\includegraphics[scale=0.25]{end_vor070501Loy.pdf}\hspace{0.5cm}
%%\includegraphics[scale=0.25]{end_vor070503QN.pdf}
%%\includegraphics[scale=0.25]{end_vor070511Loy.pdf}\hspace{0.5cm}
%%\includegraphics[scale=0.25]{end_vor070513QN.pdf}
%%\includegraphics[scale=0.25]{end_vor070521Loy.pdf}\hspace{0.5cm}
%%\includegraphics[scale=0.25]{end_vor070523QN.pdf}
%%\end{center}
%%\caption{CVT meshes with 2562 cells, generated by the Lloyd (left) and the Lloyd P-LBFGS (right), with three typical point-density functions used in the convergence efficiency tests: X3, X16, X64 (from top to bottom).}
%%\label{VTcomp}
%%\end{figure}

%%%%%%%%%%%%%%%%%%%%%%%
%%%%%%%%%%%%%%%%%%%%%%%
\subsection{Comparison of Lloyd, LBFGS and Lloyd P-LBFGS in serial}\label{serialRes}
%%%%%%%%%%%%%%%%%%%%%%%
%%%%%%%%%%%%%%%%%%%%%%%

Lloyd's method is the best-known algorithm for SCVT construction and LBFGS is a sound way to speed up the Lloyd iteration, but we believe that  the Lloyd P-LBFGS method will further speed up LBFGS, especially for highly variable multi-resolution meshes. To compare the performance of these three methods, we execute a sequence of  tests for different scenarios for serial implementations. In these tests, we only work on coarse meshes with 2562 and 10242 SCVT generators and all the runs are effected on a PC with a 2.70GHz Intel i7-3740QM CPU. The mesh points are initialized by Monte Carlo sampling based on the corresponding point-density functions. To compare both the efficiency and the mesh quality, we set the maximum number of iterations to be 2,000 and the tolerance on the size of the point movement to be $5\times 10^{-4}$. Except for that, the quasi-Newton iterations will stop when $\lVert\nabla\mathcal{F}^n\rVert / \mathcal{F}^n \leq 5\times 10^{-4}$ or $|\mathcal{F}^{n} - \mathcal{F}^{n-1}| / \mathcal{F}^{n-1} < 10^{-7}$. The tolerance $10^{-7}$ on the control of the energy decrease is to avoid wandering.  {For the integration over a triangle, we use a 4-point quadrature rule in the X3 and X16 cases, and a 9-point quadrature in the X64 case which includes extremely coarse cells.} Table \ref{tabCompSerial} provides detailed information of the three tested methods with 2,562 generators. Because the SCVT energy functional has many local minimizers, different optimization methods with even the same initial guess can still lead to different minimizers. For this reason, the final gradient norm $\lVert\nabla\mathcal{F}\rVert$ plays a role in judging the quality of final SCVT mesh.

\begin{table}[!h]
\begin{center}
\caption{Comparison of the Lloyd, LBFGS, and Lloyd P-LBFGS methods in serial with 2,562 generators.}
\label{tabCompSerial}
\begin{tabular}{llrrrll}
\toprule
Density & Method & \# iter.  & \# $\mathcal{F}$ eval. & Time (sec.) & Final $\mathcal{F}$ & Final $\lVert\nabla\mathcal{F}\rVert$\\
\midrule
\multirow{3}{*}{X3} & Lloyd & 1235 & --- & 113.869 & 1.32665e-03 & 2.6344e-05  \\
&LBFGS  &  204 & 207 & 19.965 & 1.32598e-03 & 5.4127e-07 \\
&Lloyd P-LBFGS &  170 & 171 & 16.702 & 1.32537e-03 & 1.9154e-06  \\
\midrule
\multirow{3}{*}{X16} &Lloyd  & 2000 & --- & 255.361 & 3.91468e-05 & 8.0502e-07   \\
&LBFGS  &  787 & 813 & 111.521 & 3.90003e-05 & 9.5421e-08  \\ 
&Lloyd P-LBFGS  & 309 & 334 & 45.230 & 3.89749e-05 & 5.6632e-08    \\ 
\midrule
\multirow{3}{*}{X64} & Lloyd  & 1977 & --- & 338.103 & 1.04809e-04 & 2.1095e-06   \\  
&LBFGS  & 1201 & 1290 & 236.211 & 1.08224e-04 & 1.6599e-06    \\
&Lloyd P-LBFGS  &  257 & 261 & 47.443 &1.04699e-04 & 8.6221e-08   \\
\bottomrule
\end{tabular}
\end{center}
\end{table}

\begin{figure}[!h]
\begin{center}
\begin{minipage}{0.32\textwidth}
\centering{X3 point-density}
\includegraphics[scale=0.34]{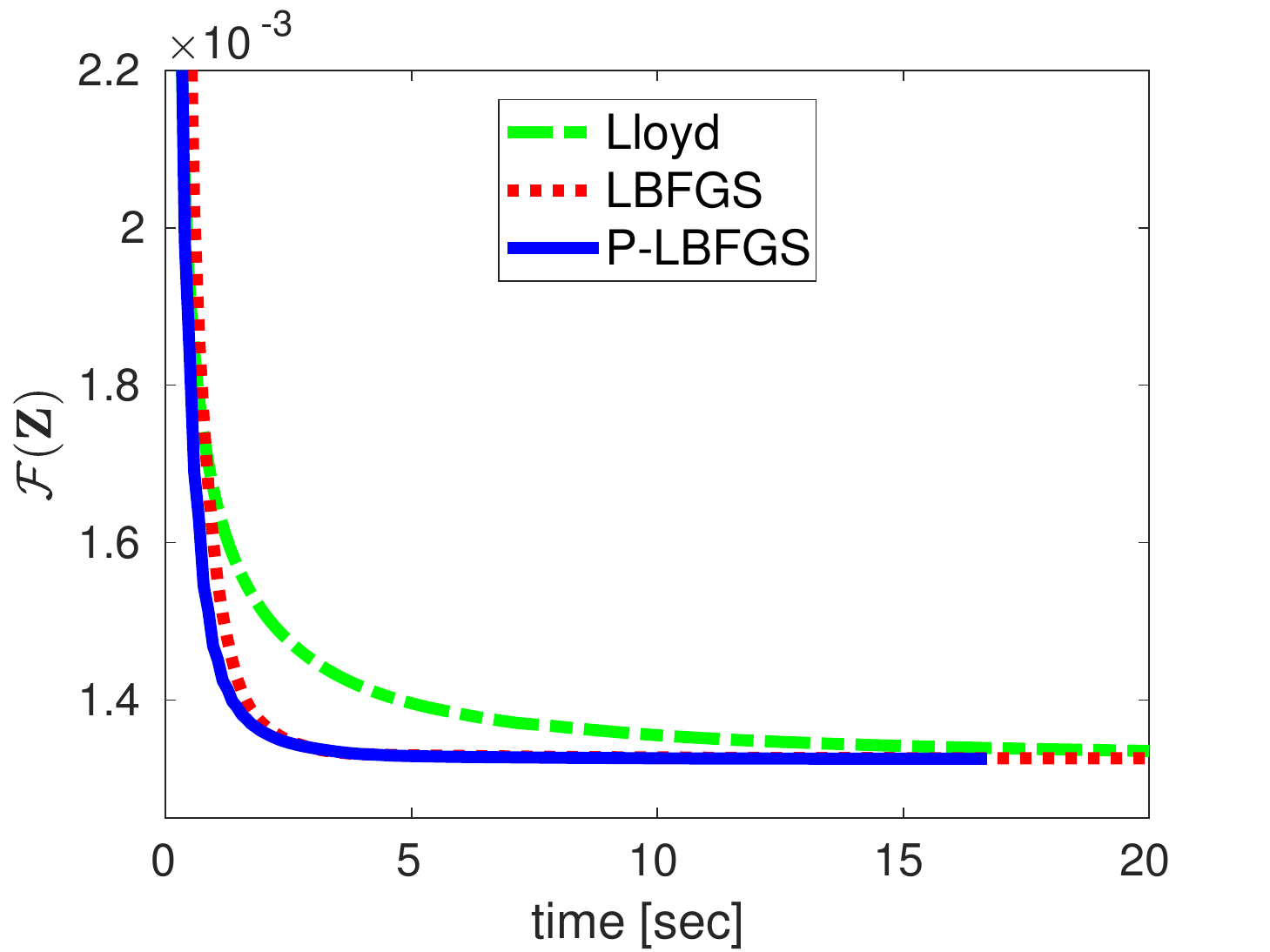} 
\end{minipage}
\begin{minipage}{0.32\textwidth}
\centering{X16 point-density}
\includegraphics[scale=0.34]{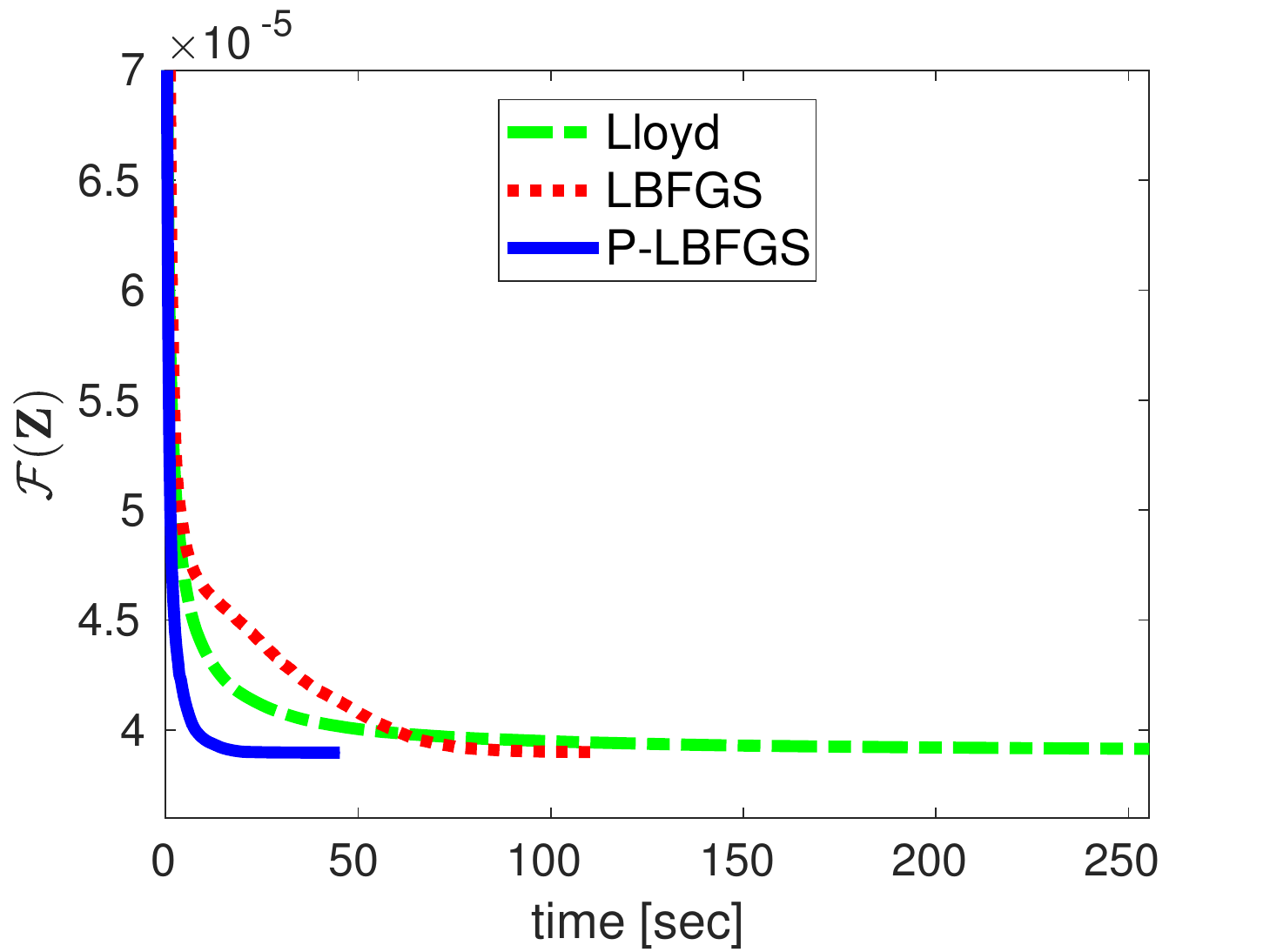} 
\end{minipage}
\begin{minipage}{0.32\textwidth}
\centering{X64 point-density}
\includegraphics[scale=0.34]{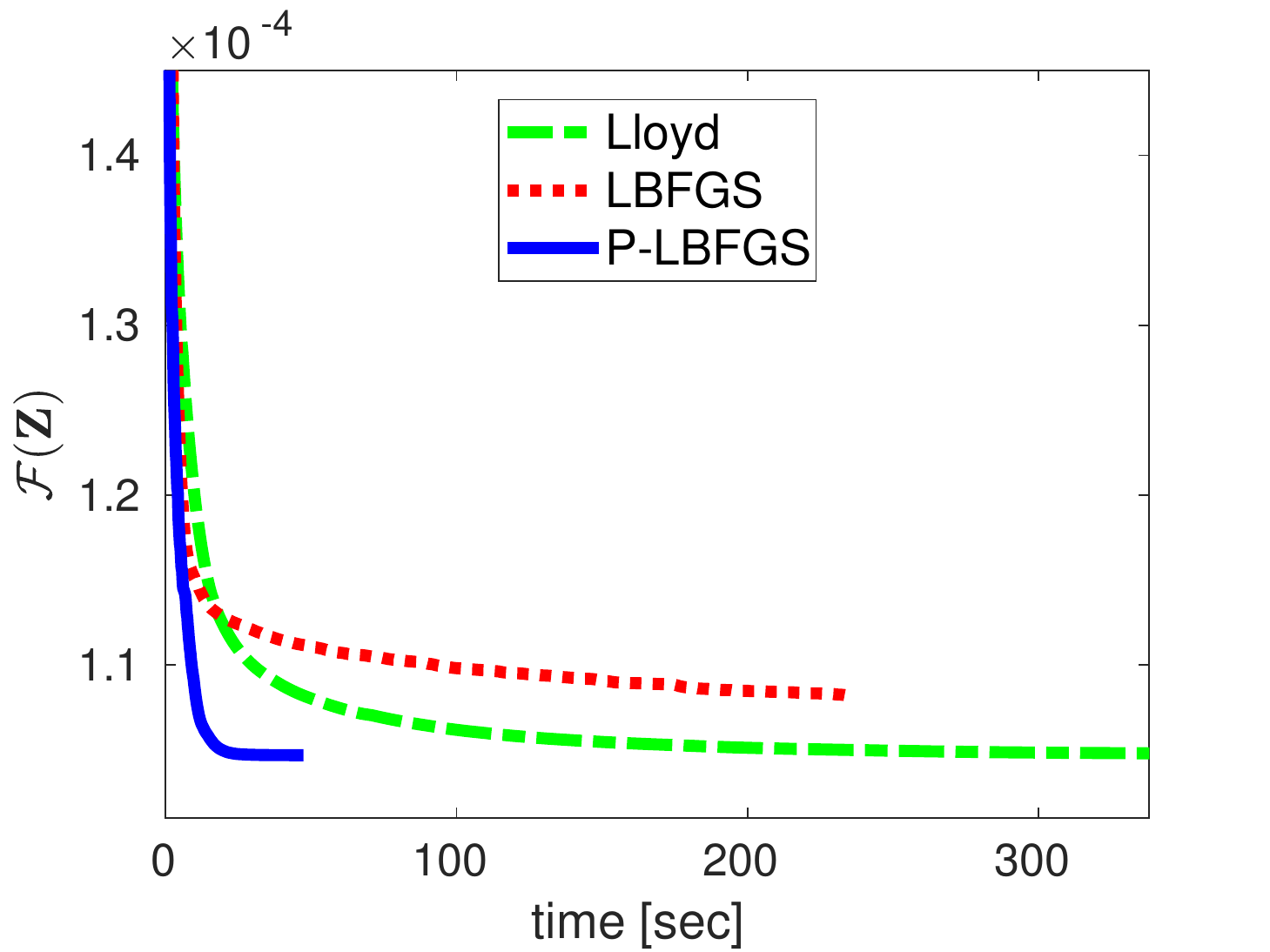} 
\end{minipage}
\begin{minipage}{0.32\textwidth}
\includegraphics[scale=0.34]{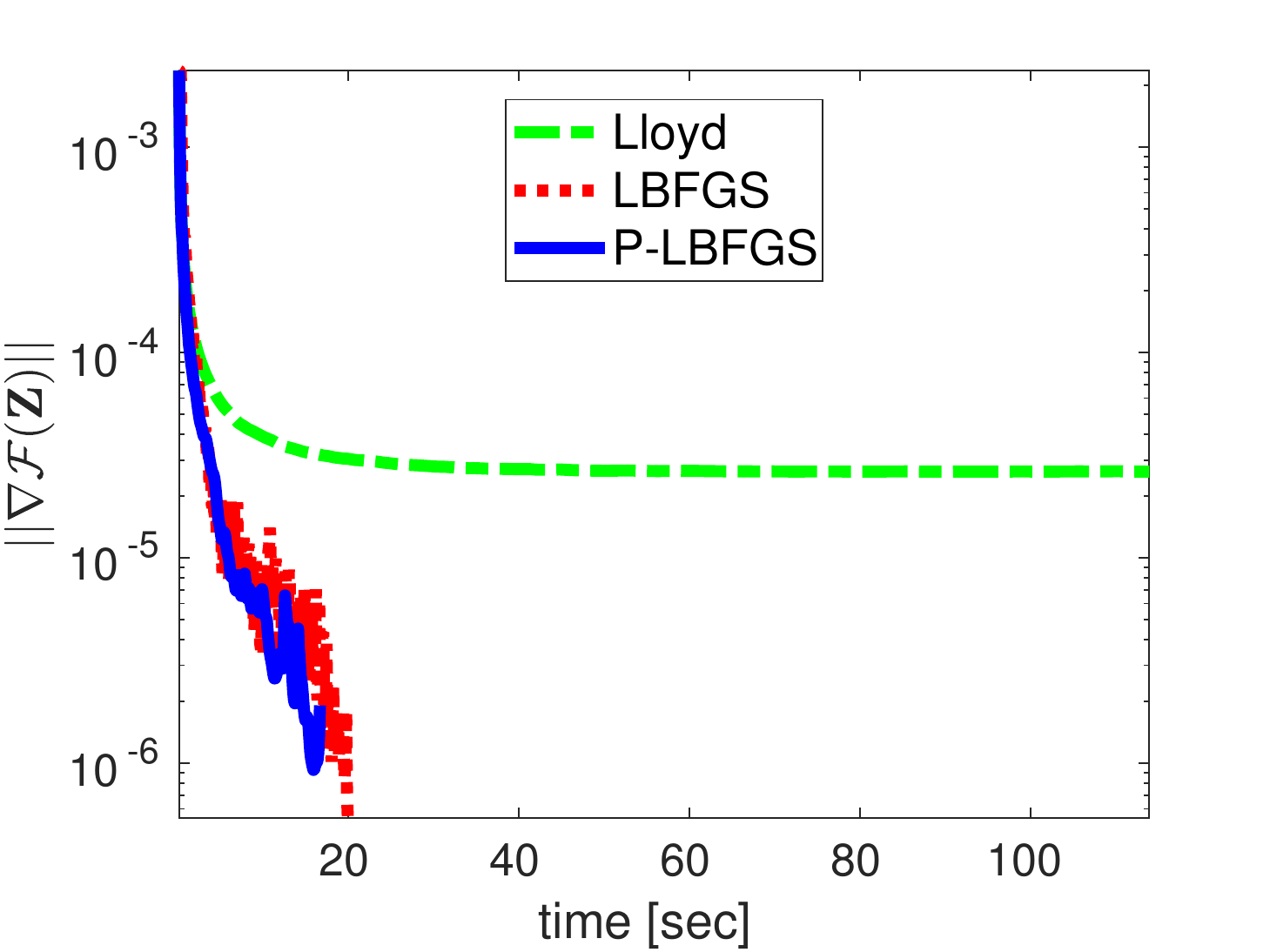} 
\end{minipage}
\begin{minipage}{0.32\textwidth}
\includegraphics[scale=0.34]{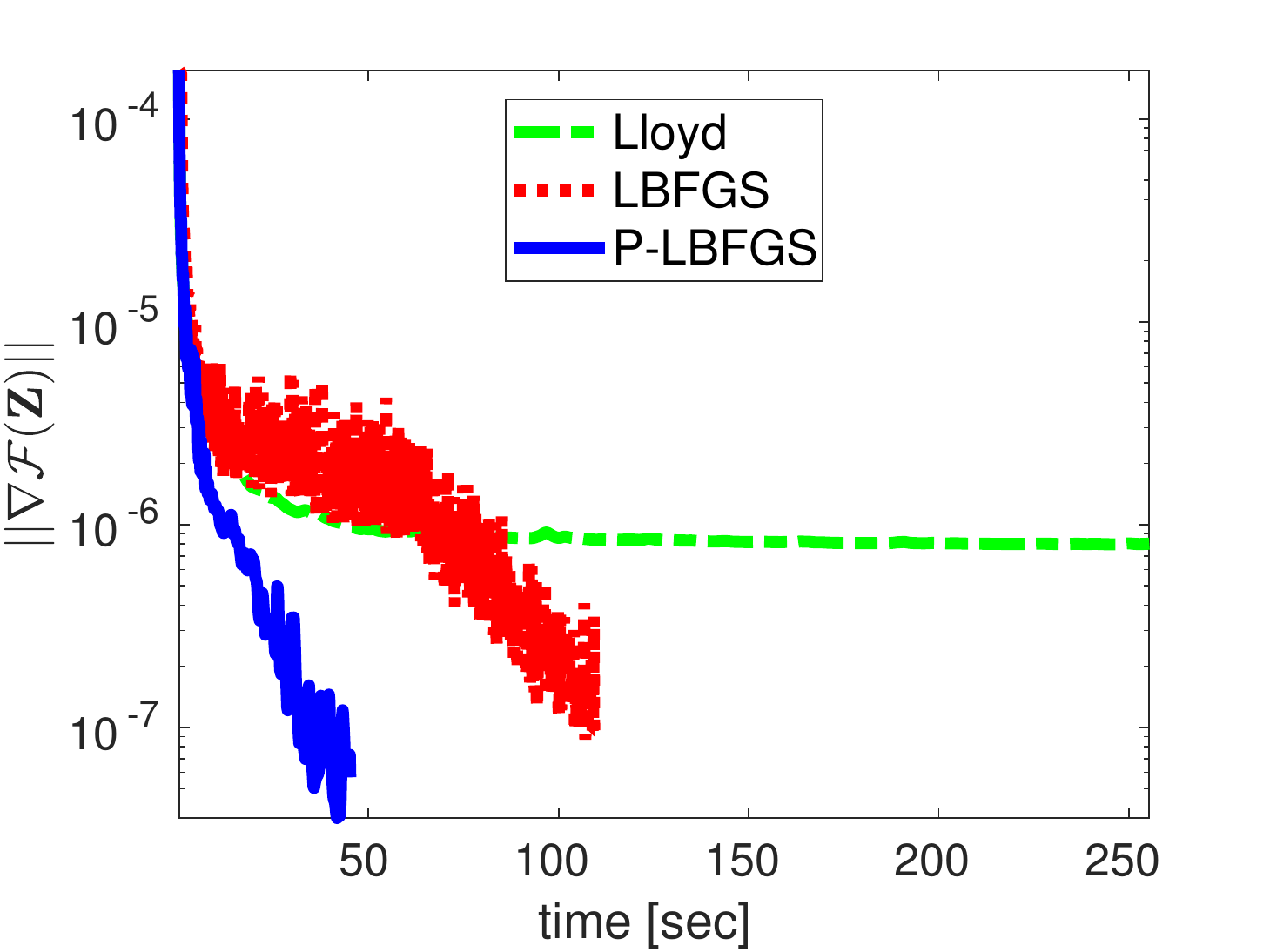} 
\end{minipage}
\begin{minipage}{0.32\textwidth}
\includegraphics[scale=0.34]{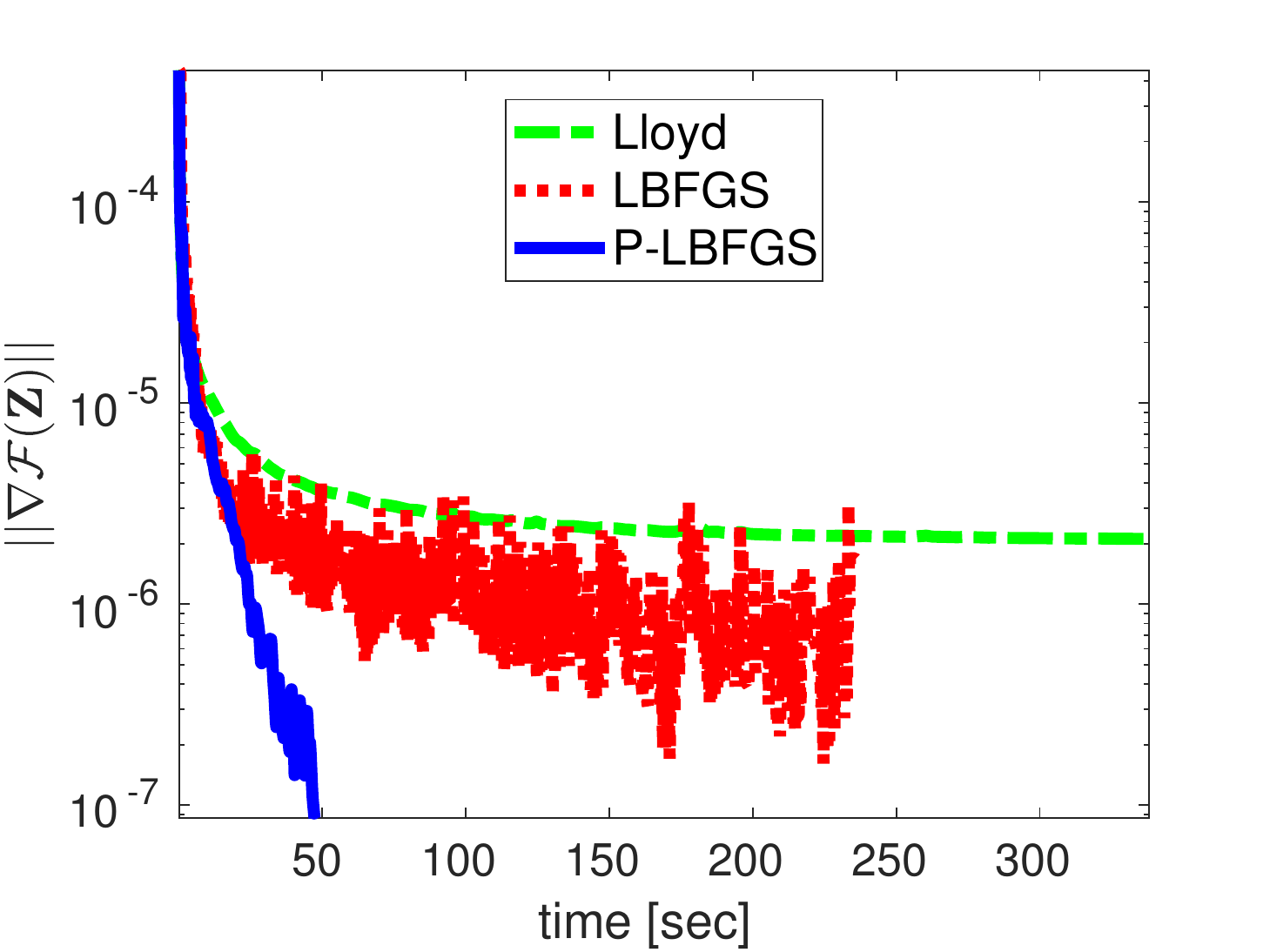} 
\end{minipage}
\begin{minipage}{0.32\textwidth}
\includegraphics[scale=0.35]{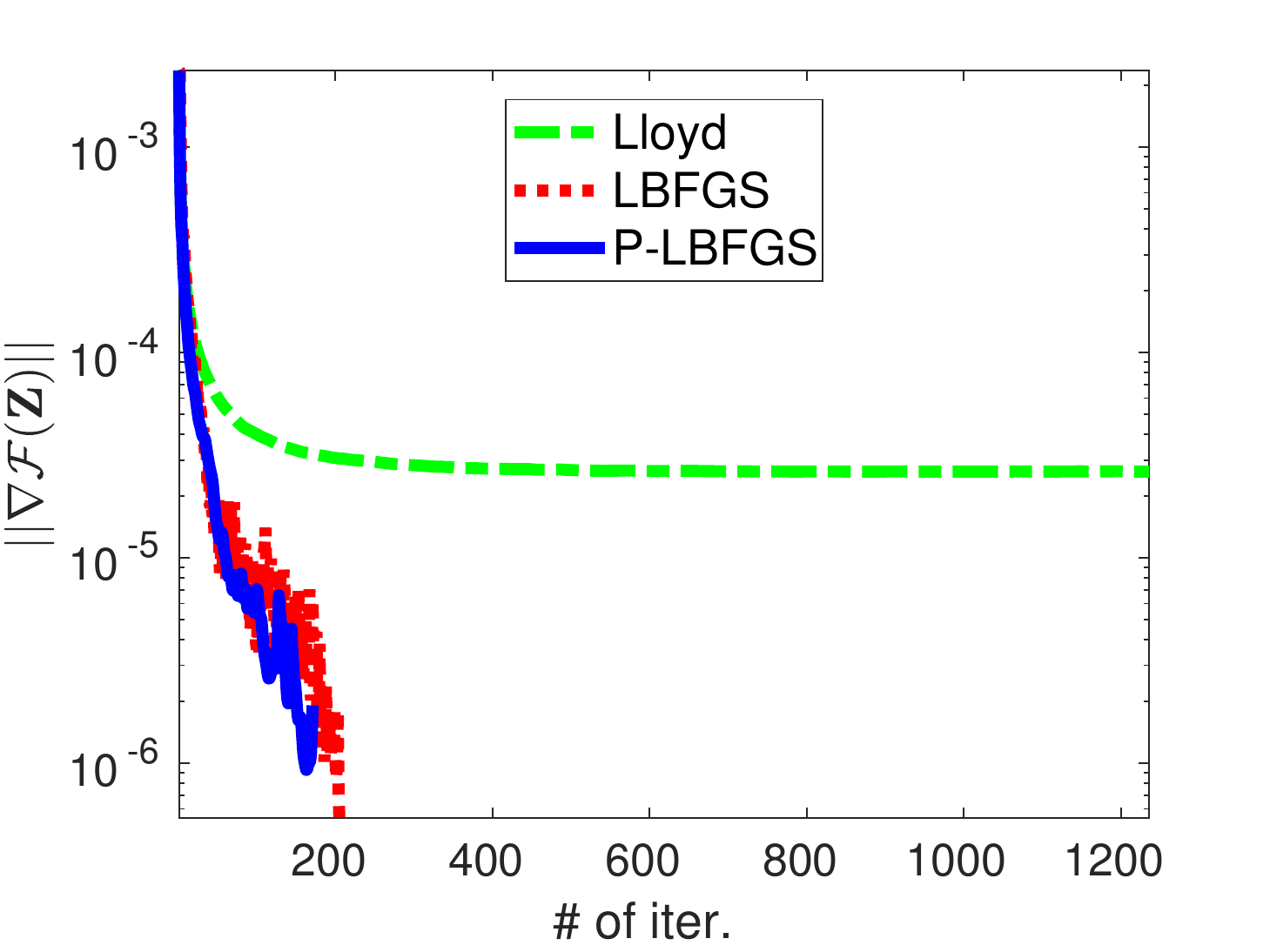} 
\end{minipage}
\begin{minipage}{0.32\textwidth}
\includegraphics[scale=0.35]{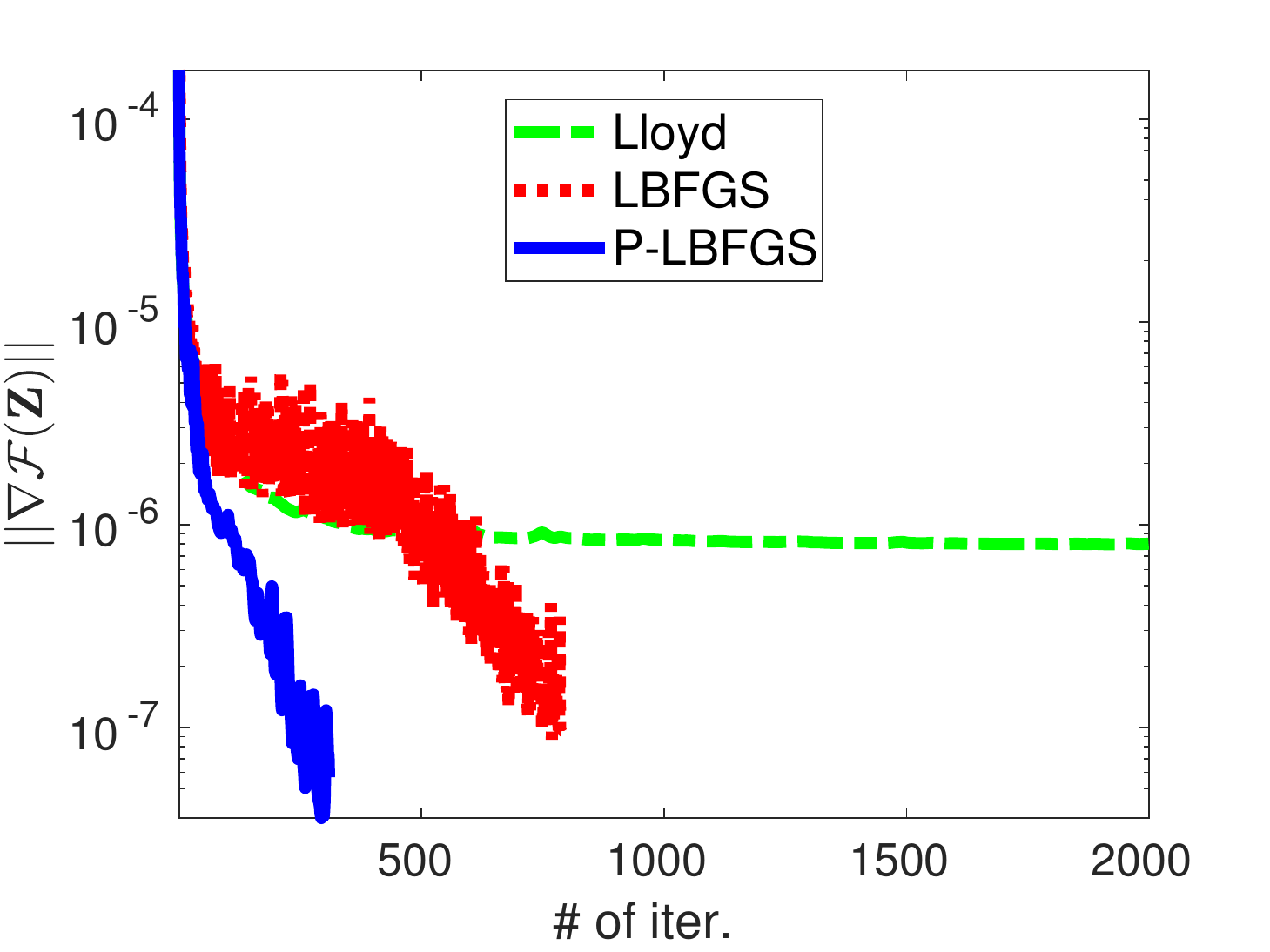} 
\end{minipage}
\begin{minipage}{0.32\textwidth}
\includegraphics[scale=0.35]{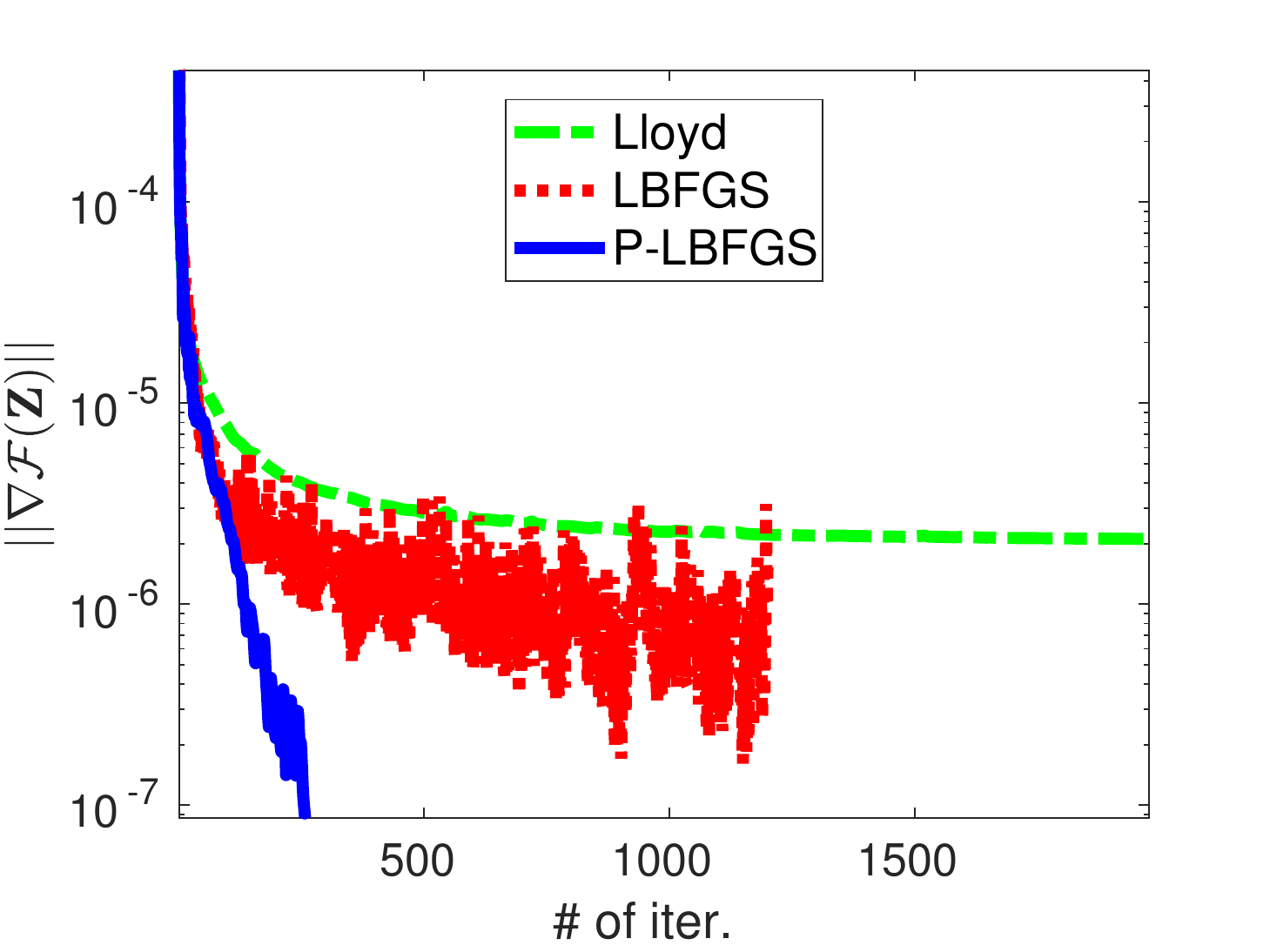} 
\end{minipage}
\caption{Plots of the performance  of the Lloyd, LBFGS, and Lloyd P-LBFGS methods in serial with 2,562 generators corresponding to Table \ref{tabCompSerial}. From top to bottom: $\mathcal{F}$ vs. computational time, $\lVert\nabla\mathcal{F}\rVert$ vs. computational time, $\lVert\nabla\mathcal{F}\rVert$ vs. \# of iterations; from left to right: X3, X16, X64 densities.}
\label{serCompPlot}
\end{center}
\end{figure}

Table \ref{tabCompSerial} shows that for all three point-density functions, Lloyd P-LBFGS is much more efficient than LBFGS and Lloyd (1.18 to 5.02 times faster than LBFGS and 5.66 to 7.19 times faster than Lloyd). We also observed that the relative efficiency of Lloyd P-LBFGS becomes even better in the case with 10,242 generators, as illustrated in Table \ref{tabCompSerial10242}. Among the three methods Lloyd P-LBFGS always returns the smallest value on the energy function and nearly so for the gradient norm, which usually mean better mesh quality. The speed-up of LBFGS by the Lloyd preconditioner is more obvious as we increase the mesh size ratio; this is also shown by the plots of the iteration histories in Figure~\ref{serCompPlot} which provides the numerical results in terms of the energy functional and the gradient norm with respect to the number of iterations and the computation time. 

\begin{table}[!h]
\begin{center}
\caption{Comparison of the Lloyd, LBFGS, and Lloyd P-LBFGS methods in serial with 10,242 generators.}
\label{tabCompSerial10242}
\begin{tabular}{llrrrll}
\toprule
Density & Method & \# iter.  & \# $\mathcal{F}$ eval. & Time (sec.) & Final $\mathcal{F}$ & Final $\lVert\nabla\mathcal{F}\rVert$\\
\midrule
\multirow{3}{*}{X3} & Lloyd & 2000 & --- & 851.965 & 3.31424e-04 & 3.3334e-06  \\ % 0705011Loy
&LBFGS  &  347 & 354 & 145.716 & 3.31306e-04 & 1.6461e-07 \\ % 0705022QN
&Lloyd P-LBFGS &  315 & 316 & 130.262 & 3.31325e-04 & 1.5729e-07  \\ % 0705033QN
\midrule
\multirow{3}{*}{X16} &Lloyd  & 2000 & --- & 1138.490 & 9.93638e-06 & 1.1110e-07   \\ % 0705111Loy
&LBFGS  &  1608 & 1621 & 1089.000 & 9.76267e-06 & 2.8646e-08  \\ % 0705122QN
&Lloyd P-LBFGS  & 533 & 534 & 335.890 & 9.75698e-06 & 4.8184e-09    \\ % 0705133QN
\midrule
\multirow{3}{*}{X64} & Lloyd  & 2000 & --- & 1406.740 & 2.65790e-05 & 3.6220e-07   \\  % 1028211Loy
&LBFGS  & 2000 & 2024 & 1504.600 & 2.73630e-05 & 1.1220e-07    \\ % 1028222QN
&Lloyd P-LBFGS  &  405 & 417 & 320.680 & 2.61190e-05 & 1.7556e-08   \\ % 1028233QN
\bottomrule
\end{tabular}
\end{center}  
\end{table}

% figures generated by run_report_compareIter_serial3rd.m

\begin{figure}[!h]
\begin{center}
\includegraphics[scale=0.5]{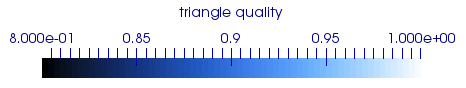}\\
\includegraphics[scale=0.35]{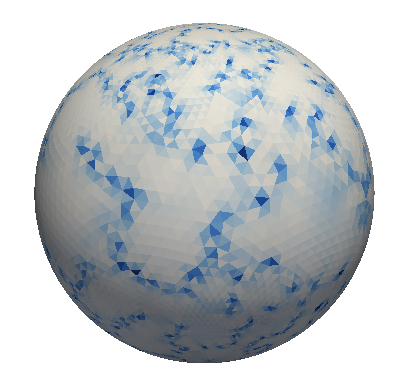}
\includegraphics[scale=0.35]{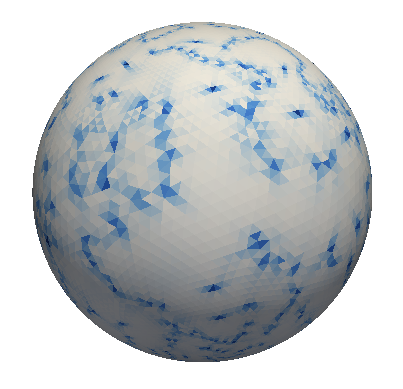}
\includegraphics[scale=0.35]{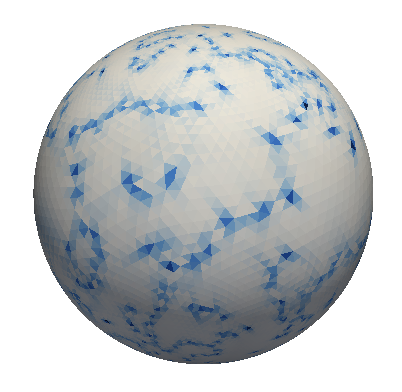}

\includegraphics[scale=0.35]{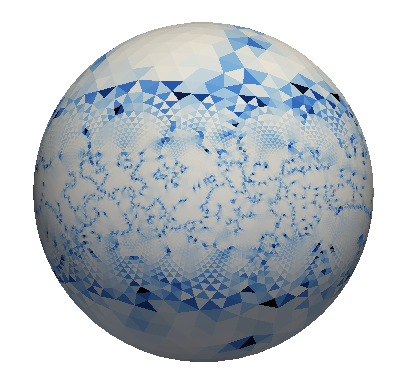}
\includegraphics[scale=0.35]{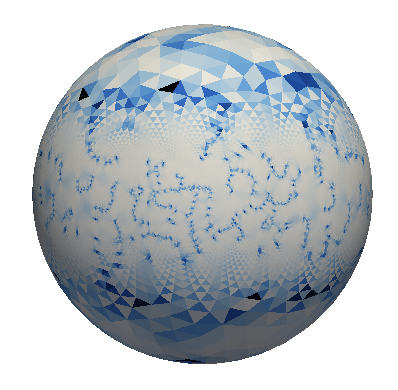}
\includegraphics[scale=0.35]{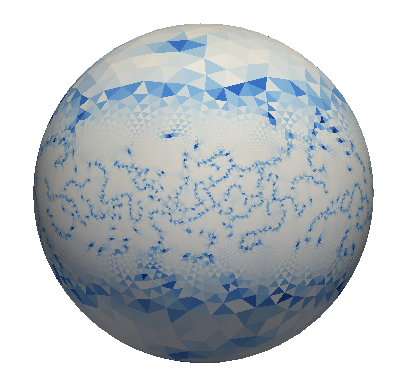}
\end{center}
\caption{Distribution of triangle qualities of the produced SCVT meshes with 10,242 generators. From top to bottom: X3, X16 densities; from left to right: Lloyd, LBFGS, Lloyd P-LBFGS.}
\label{TriQserial}
\end{figure}

SCVT meshes with 2,562 generators generated by Lloyd P-LBFGS are plotted in Figure~\ref{VTcomp} from which we can see that the mesh points are very regularly distributed and that the mesh transitions from coarse to fine regions are very smooth. The quality of SCVT meshes can also be characterized by the triangle quality measure ranging from 0 to 1; see the definition in (\ref{Qdef}). In Figure~\ref{TriQserial}, we plot the distribution of the triangle  quality of the SCVT meshes with 10,242 generators produced by the three methods; we here only plot the case of X3 and X16 with 10,242 generators because the X64 case has some very triangular cells that are too coarse to be properly visualized by the Paraview software. In this figure, the lighter the color, the better the quality. The results show that in general both LBFGS and Lloyd P-LBFGS produce better SCVT meshes in term of quality than Lloyd and for the X16 case, the Lloyd P-LBFGS clearly provides meshes with better quality than other two methods.

%%%%%%%%%%%%%%%%%%%%%%%
%%%%%%%%%%%%%%%%%%%%%%%
\subsection{Comparison of Lloyd  and graph Laplacian preconditioners in P-LBFGS in serial}\label{GLvsLP}
%%%%%%%%%%%%%%%%%%%%%%%
%%%%%%%%%%%%%%%%%%%%%%%

We next make a direct performance comparison of the Lloyd preconditioner  with the graph Laplacian preconditioner in the P-LBFGS scheme.
The graph Laplacian is computed with a Cholesky decomposition solver as done in \cite{Hateley2015}.
Because the Laplacian matrix is singular, the numerical solution is unique up to a constant; depending on which graph Laplacian preconditioner is used, the corresponding P-LBFGS method could perform very differently. Some solutions of the Laplacian system may lead to a direction with a poor length for line search or even a non-decent direction. As a result, we will handle the singularity of the Laplacian operators by a diagonal perturbation. Two approaches are tested . One is the diagonal shifting by adding scaled identity: $\mathbf{A} \leftarrow\mathbf{A}+10^{-6}\mathbf{I}$. We call this approach Laplacian(a6). The other one is to multiply the diagonal entries by a constant factor: $\mathbf{A}_{ii} \leftarrow  \mathbf{A}_{ii}\cdot(1+10^{-2})$. This method is denoted by Laplacian(m2). We apply the graph Laplacian preconditioner in each iteration so that the numerical results can represent, in principle, the best performance of Laplacian preconditioner. 

It is worth mentioning that the point-density function in \cite[Example 4]{Hateley2015} has min-max ratio $e^{-80}$, thus the mesh size ratio is as huge as $e^{20}$; similarly in \cite[Example 5]{Hateley2015}, the mesh size ratio is $e^{2.5}\approx12.18$. Our testing densities X16 and X64 are in a category comparable to \cite{Hateley2015}. In the following tests, the parameter configuration is exactly the same as in Section~\ref{serialRes}. We investigate, in different scenarios, not only the necessary number of iterations for  convergence but also the final mesh quality which is an even higher priority in many applications \cite{ringler2011exploring, Sakaguchi2015}. 

\begin{figure}[!h]
\begin{center}
\begin{minipage}{0.32\textwidth}
\centering{X3 point-density}
\includegraphics[scale=0.35]{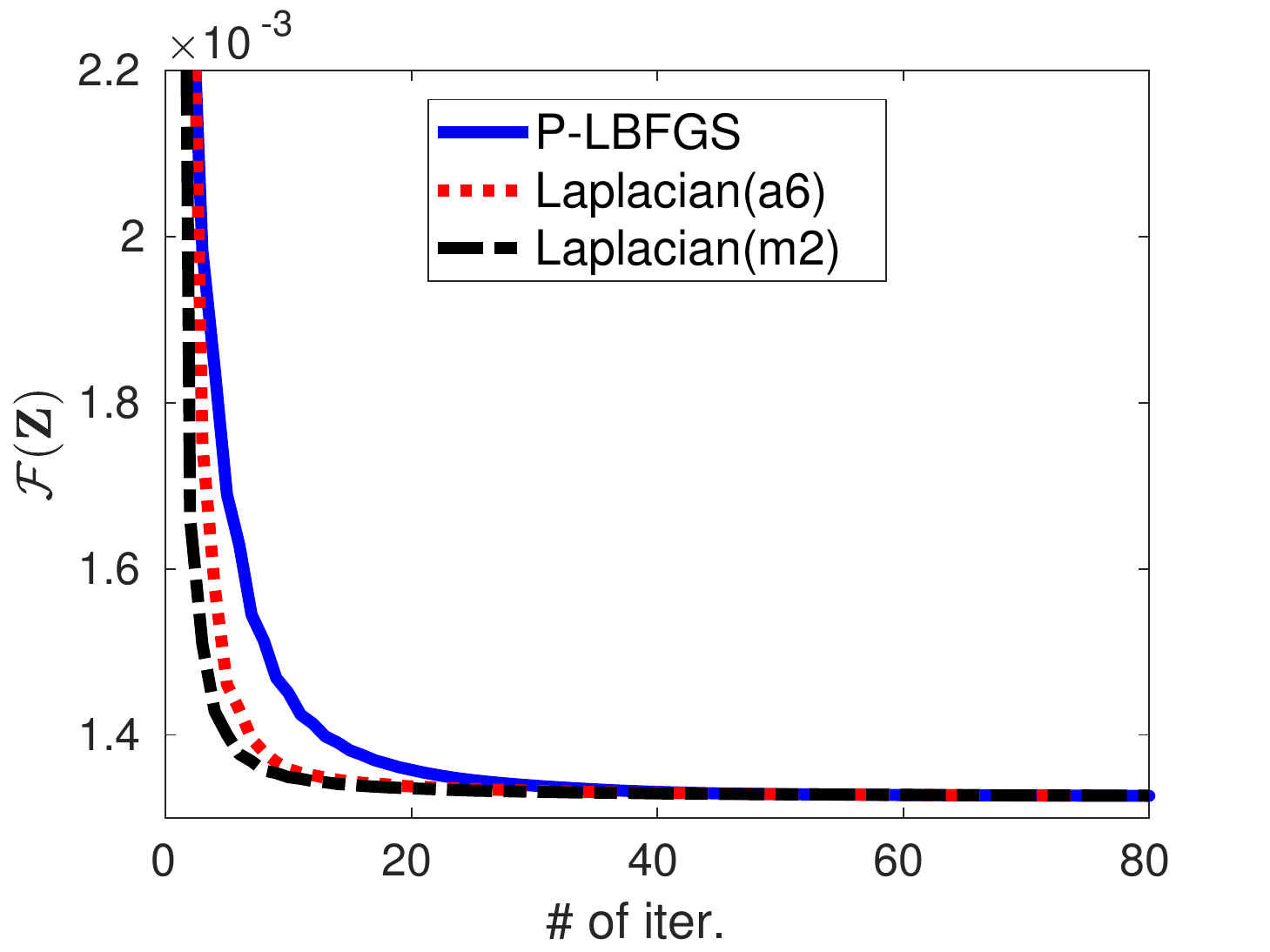} 
\end{minipage}
\begin{minipage}{0.32\textwidth}
\centering{X16 point-density}
\includegraphics[scale=0.35]{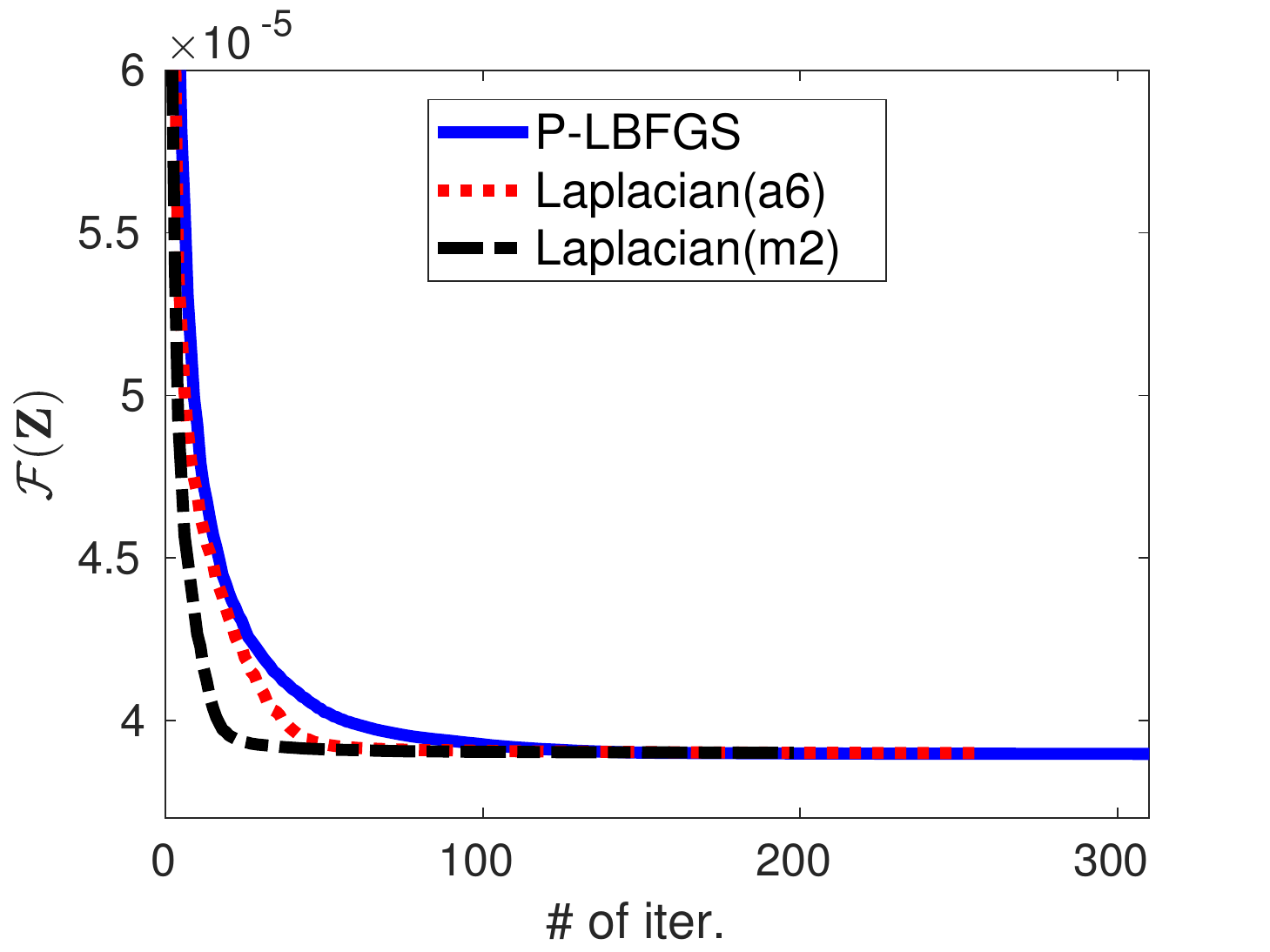} 
\end{minipage}
\begin{minipage}{0.32\textwidth}
\centering{X64 point-density}
\includegraphics[scale=0.35]{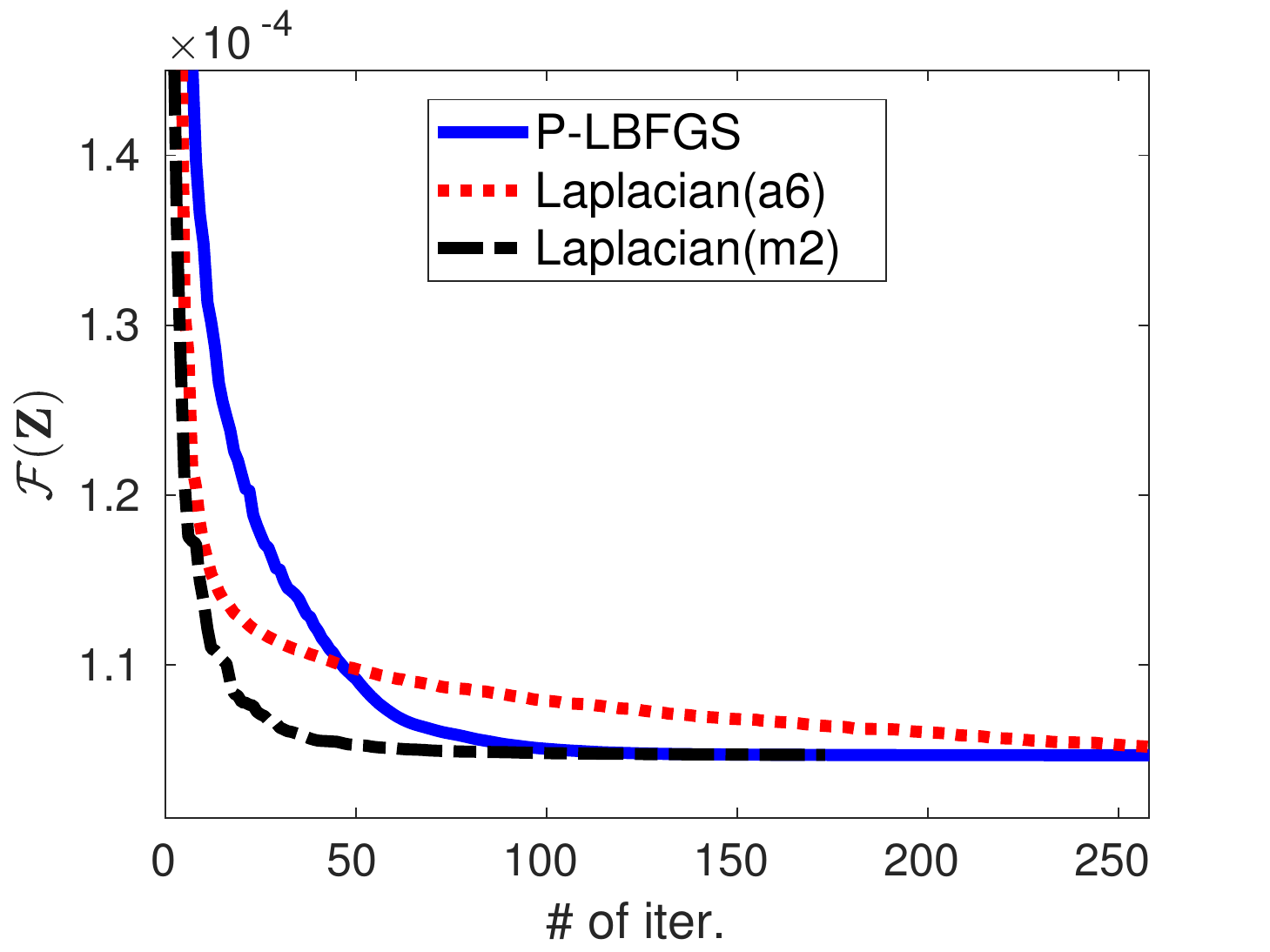} 
\end{minipage}
\begin{minipage}{0.32\textwidth}
\includegraphics[scale=0.35]{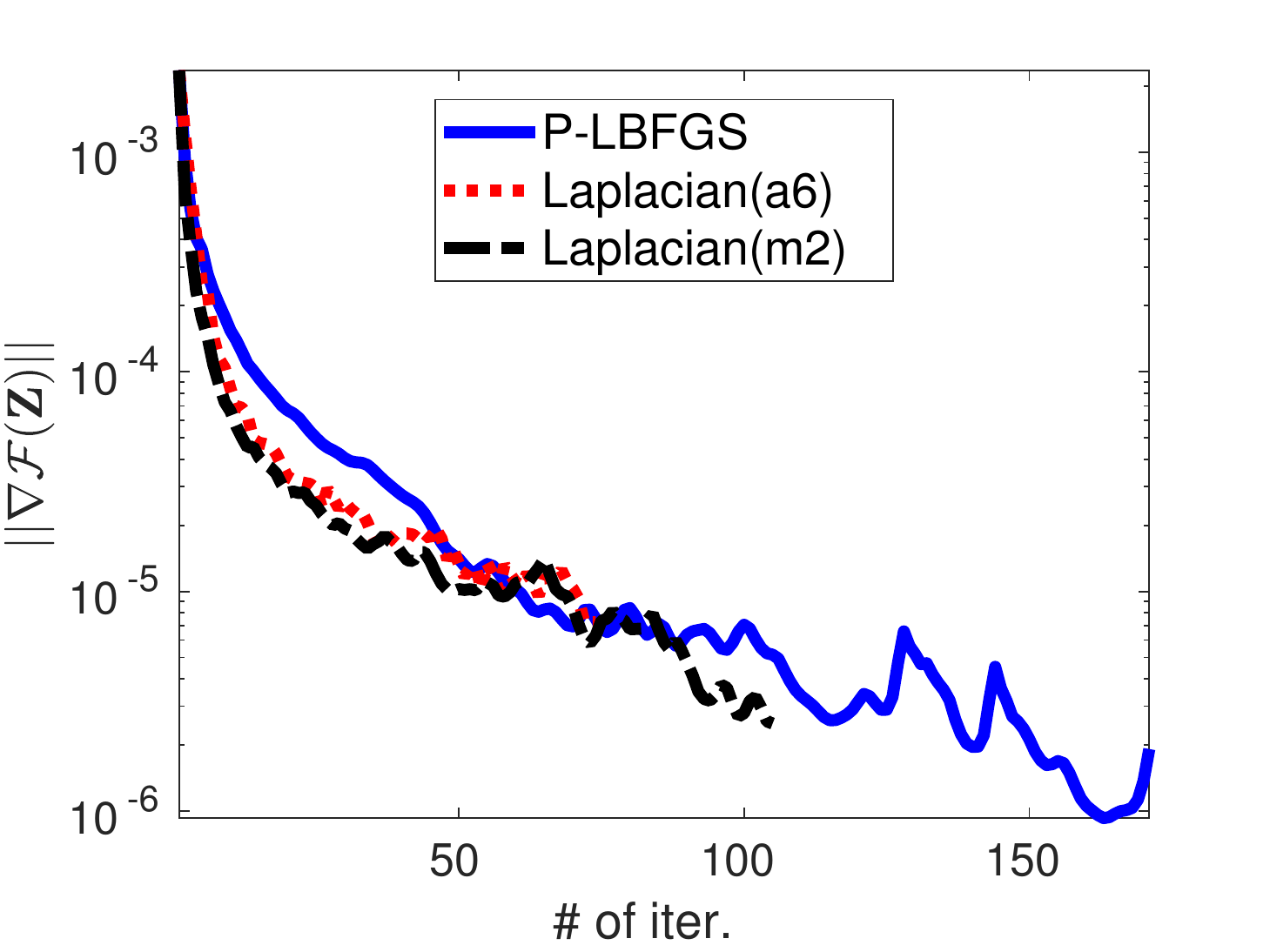} 
\end{minipage}
\begin{minipage}{0.32\textwidth}
\includegraphics[scale=0.35]{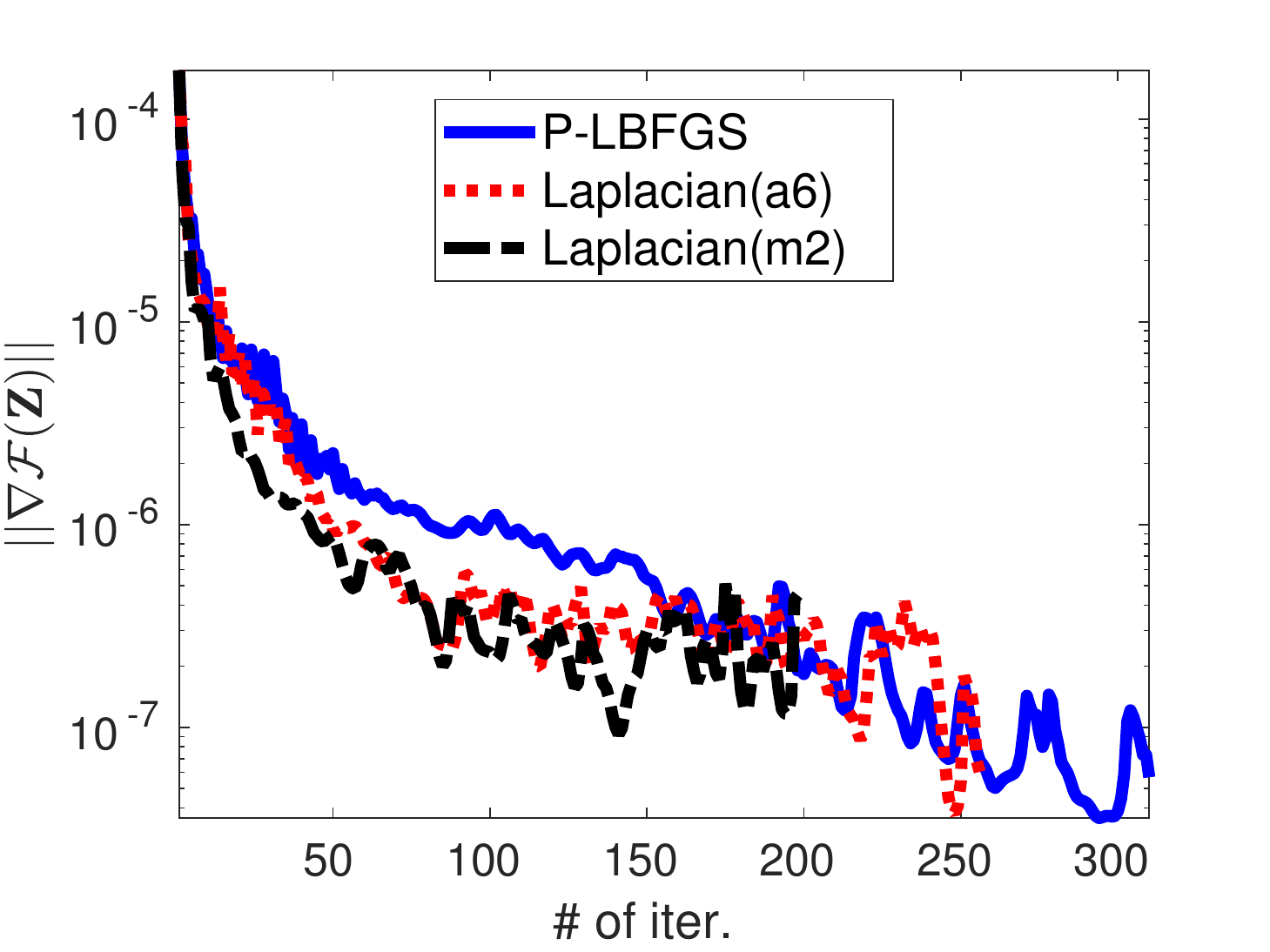} 
\end{minipage}
\begin{minipage}{0.32\textwidth}
\includegraphics[scale=0.35]{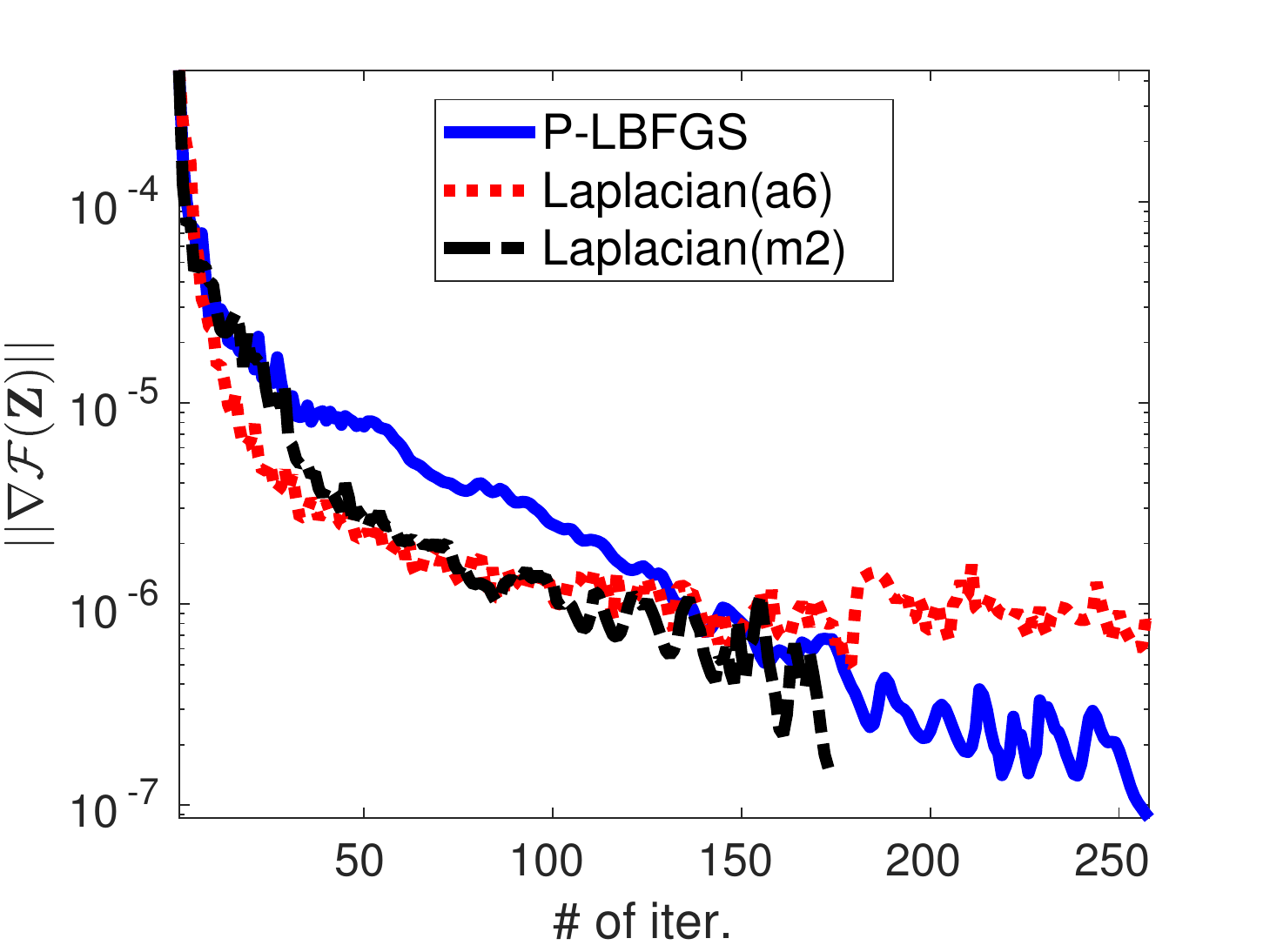} 
\end{minipage}
%\rule[0.5ex]{\linewidth}{1pt}
\begin{minipage}{0.32\textwidth}
\includegraphics[scale=0.35]{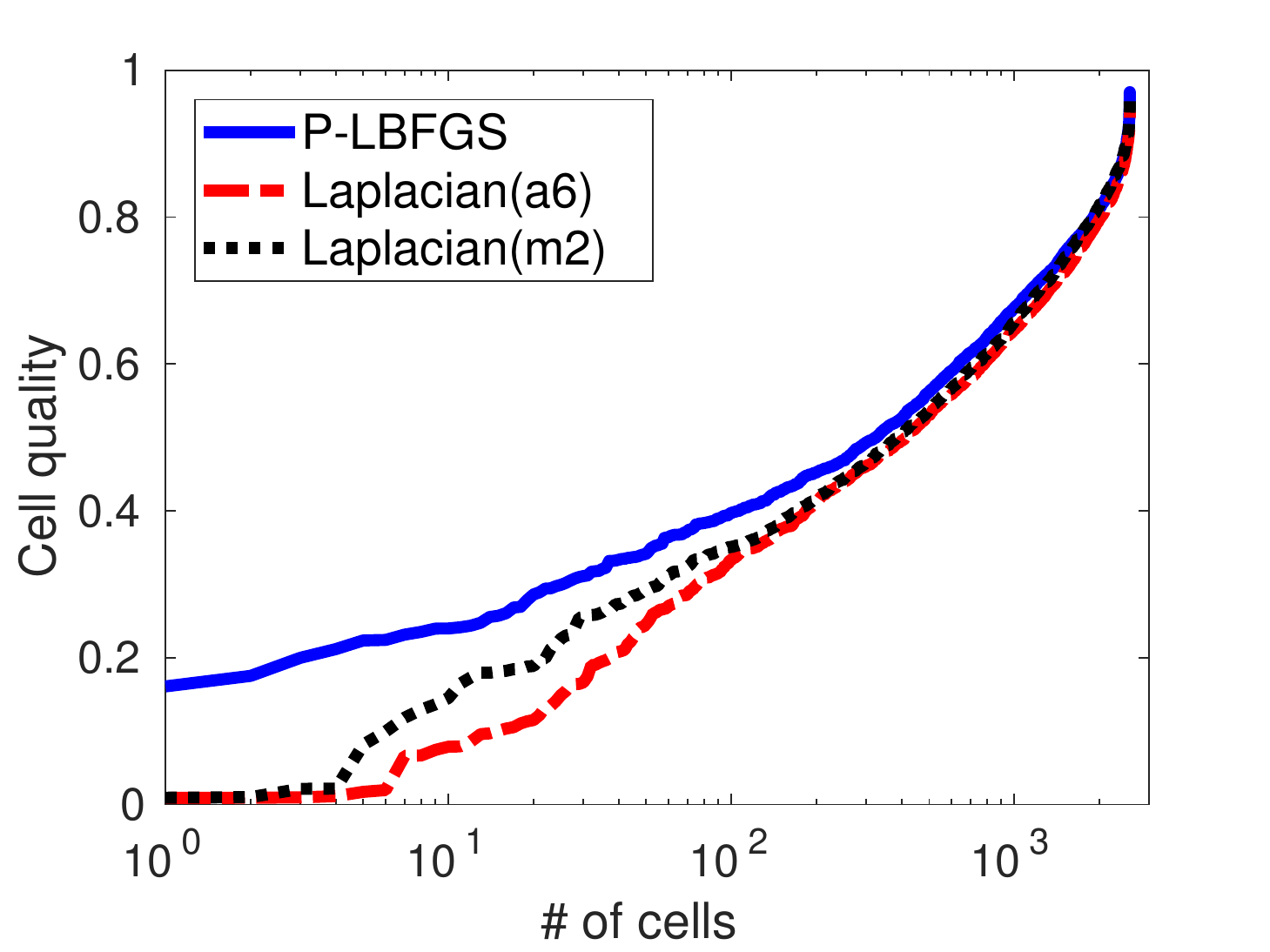} 
\end{minipage}
\begin{minipage}{0.32\textwidth}
\includegraphics[scale=0.35]{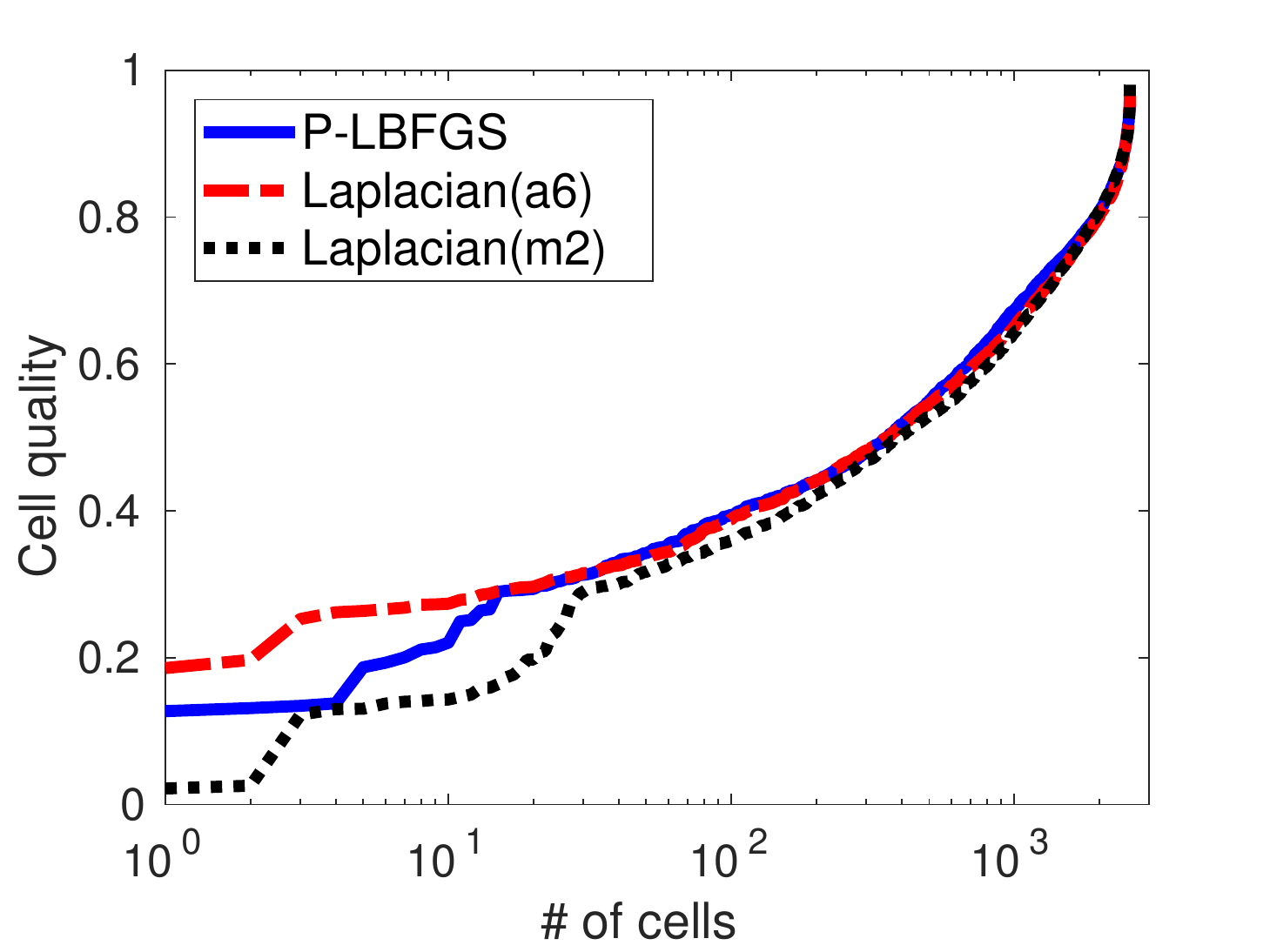} 
\end{minipage}
\begin{minipage}{0.32\textwidth}
\includegraphics[scale=0.35]{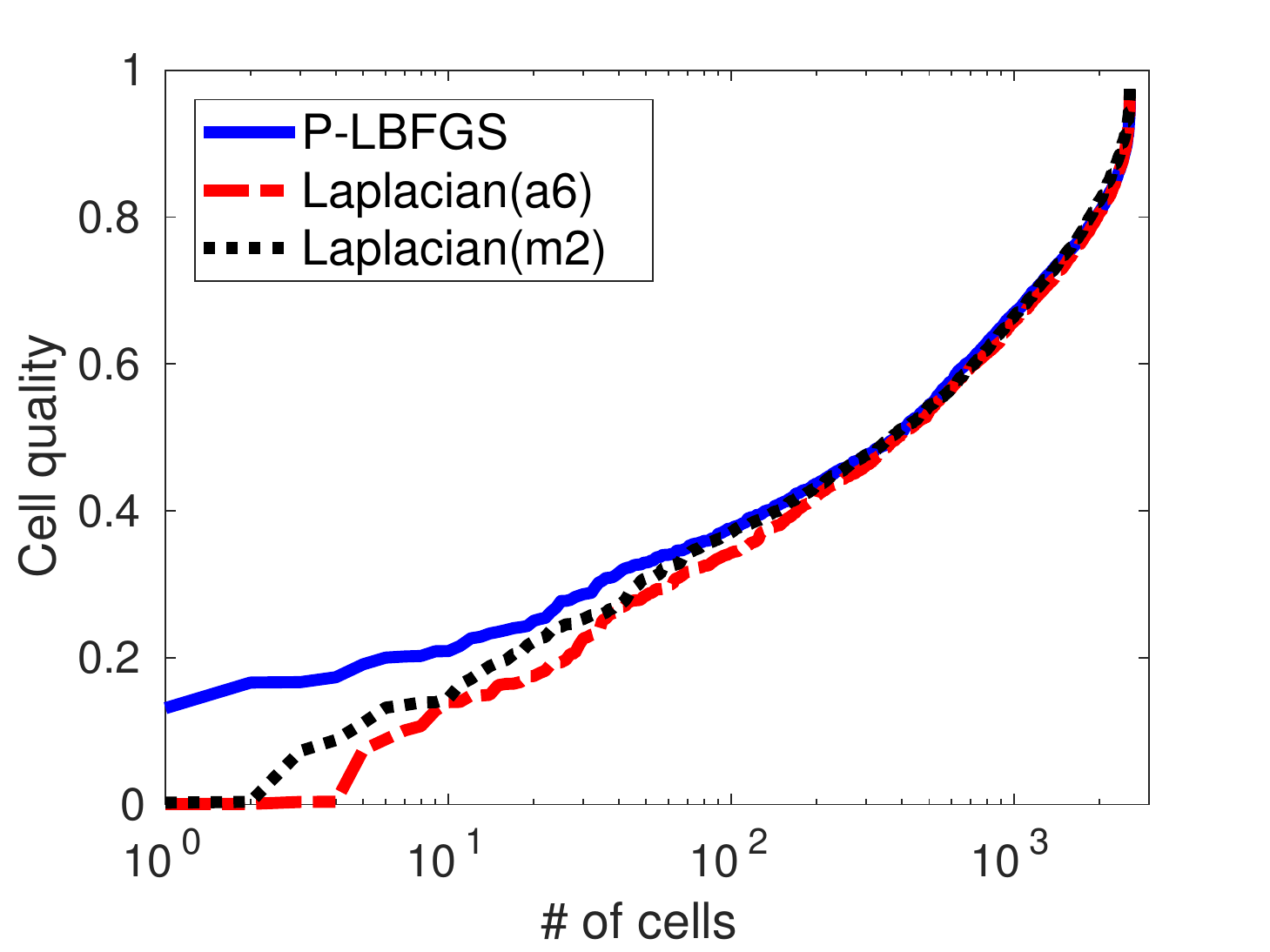} 
\end{minipage}
\caption{Plots of the performance of the Lloyd preconditioner (denoted by P-LBFGS) and the graph Laplacian preconditioner with 2,562 generators, corresponding to Table 
\ref{tabCompGL}. From top to bottom: $\mathcal{F}$ vs. \# of iterations, $\lVert\nabla\mathcal{F}\rVert$ vs. \# of iterations, cell quality $cellQ$; from left to right: X3, X16, X64 
densities.}
\label{figCompGL}
\end{center}
\end{figure}

Detailed test results are provided in Table \ref{tabCompGL} and Figure~\ref{figCompGL}. We see that Laplacian(m2) has the overall fastest convergence rate if we only focus on the number of iterations without mentioning the computational time. Laplacian(a6) requires fewer iteration than the Lloyd preconditioner in the X3 test case, whereas not significantly fewer in the X16 case, and, in the X64 test case, it requires more iterations.  Thus, we take Laplacian(m2) as a better representative for the graph Laplacian preconditioners proposed in \cite{Hateley2015}. The graph Laplacian preconditioner can dramatically reduce the necessary number of iterations in P-LBFGS; this is their main advantage as it was proposed. A natural concern 
is the computational cost incurred by having to solve discrete Laplacian systems. We observe that with respect to computational time, the Lloyd preconditioner costs much less time in all test cases. 
Another aspect which we care even more about is mesh quality. As one can see from Table \ref{tabCompGL} and  Figure~\ref{figCompGL}, Laplacian(m2) always returns larger values than Lloyd P-LBFGS for the final energy  and 
the gradient norm. {We also find that in this set} of tests Laplacian(m2) actually always stops based on the criterion of $|\mathcal{F}^{n} - \mathcal{F}^{n-1}| / \mathcal{F}^{n-1}$ rather than on the gradient criterion.
It seems the final energies only have little differences, but the mesh qualities are significantly different, as illustrated by the plot of cell quality in the third row of Figure~\ref{figCompGL}. 
We conclude that Lloyd preconditioner outperforms the graph Laplacian(m2) at providing a SCVT grid with ``optimal'' quality and computational times.

\begin{table}[!h]
\begin{center}
\caption{Comparison of the Lloyd and graph Laplacian preconditioners for P-LBFGS with 2,562 generators. 
%The time count on Lloyd's method doesn't include the time for $\mathcal{F}$ and $\lVert\nabla\mathcal{F}\rVert$ computation. 
}\label{tabCompGL}
\begin{tabular}{llrrrll}
\toprule
Density  & Method  & \# iter.  & \# $\mathcal{F}$ eval. & Time(sec) & Final $\mathcal{F}$ & Final $\lVert\nabla\mathcal{F}\rVert$\\
\midrule
\multirow{3}{*}{X3}  
& Laplacian(a6) P-LBFGS& 73 & 93 & 29.503 & 1.32712e-03 & 7.2298e-06  \\
& Laplacian(m2) P-LBFGS& 104 & 106 & 32.805 & 1.32664e-03 & 2.5014e-06  \\
&Lloyd P-LBFGS &  170 & 171 & 16.701 & 1.13253e-03 & 1.9154e-06  \\
\midrule
 \multirow{3}{*}{X16} &  
 Laplacian(a6) P-LBFGS & 255 & 262 & 90.697& 3.89986e-05 & 6.0327e-08\\  
  & Laplacian(m2) P-LBFGS& 197 & 205 & 70.725 & 3.90087e-05 & 4.4906e-07 \\  
&Lloyd P-LBFGS  & 309 & 334 & 45.229 & 3.89749e-05 & 5.6632e-08    \\ 
 \midrule
 \multirow{3}{*}{X64} & Laplacian(a6) P-LBFGS& 329 & 463 &179.691 & 1.04747e-04 & 1.8575e-07 \\  
 & Laplacian(m2) P-LBFGS& 172 & 202 & 78.739 & 1.04720e-04 & 1.5299e-07  \\ 
&Lloyd P-LBFGS  &  257 & 261 & 47.443 &1.04699e-04 & 8.6221e-08   \\
\bottomrule
\end{tabular}
\end{center}
\end{table}

%%%%%%%%%%%%%%%%%%%%%%%
%%%%%%%%%%%%%%%%%%%%%%%
\subsection{Parallel scalability of the Lloyd P-LBFGS}
%%%%%%%%%%%%%%%%%%%%%%%
%%%%%%%%%%%%%%%%%%%%%%%

Now we study the scaling performance of the parallel P-LBFGS method.
All the parallel tests were run on a high performance computing (HPC) cluster with 16 Intel Xeon ``model E5-2670'' cores per node. 

\begin{figure}[h!]
\begin{center}
\includegraphics[scale=0.36]{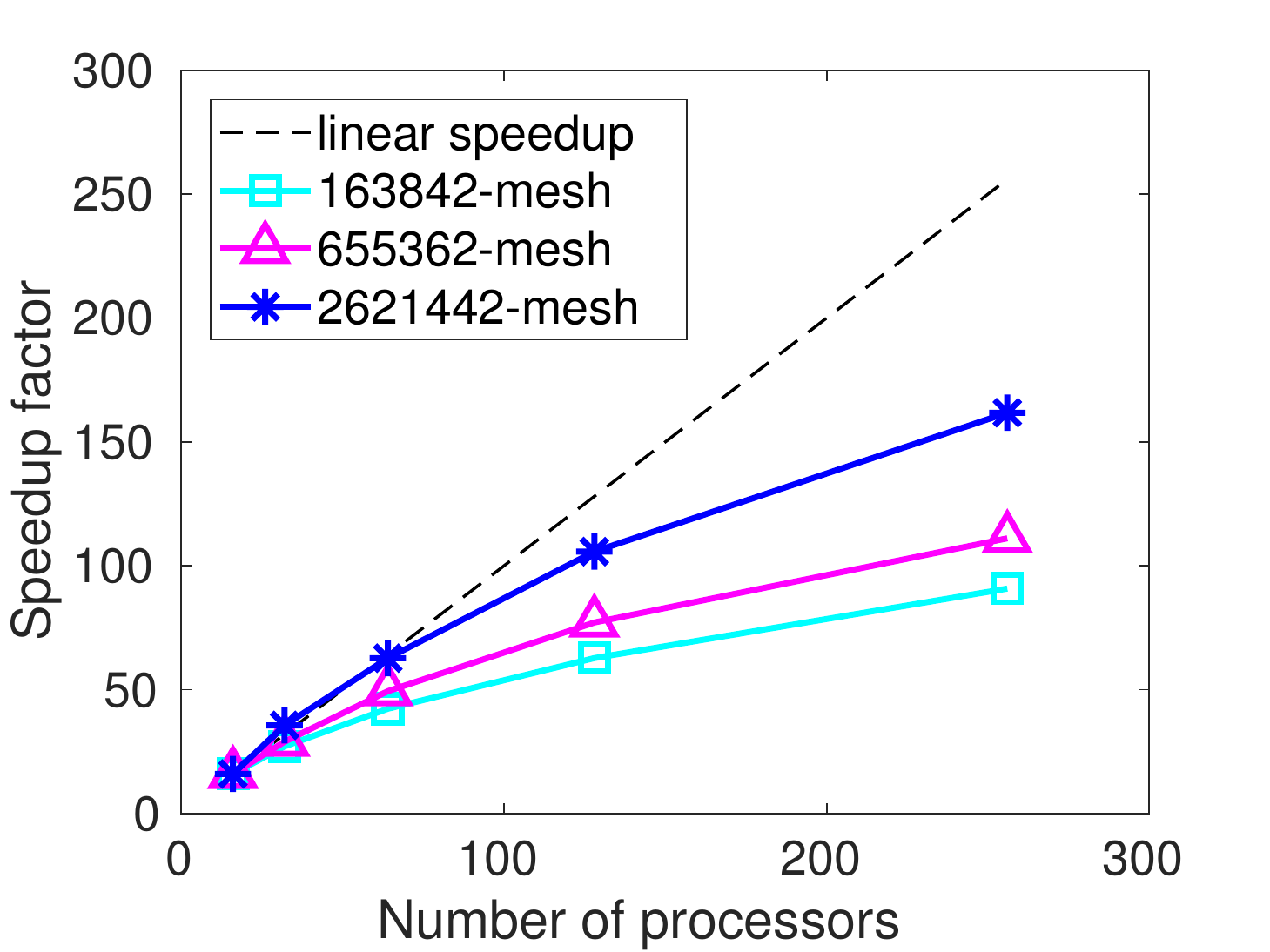}
\includegraphics[scale=0.36]{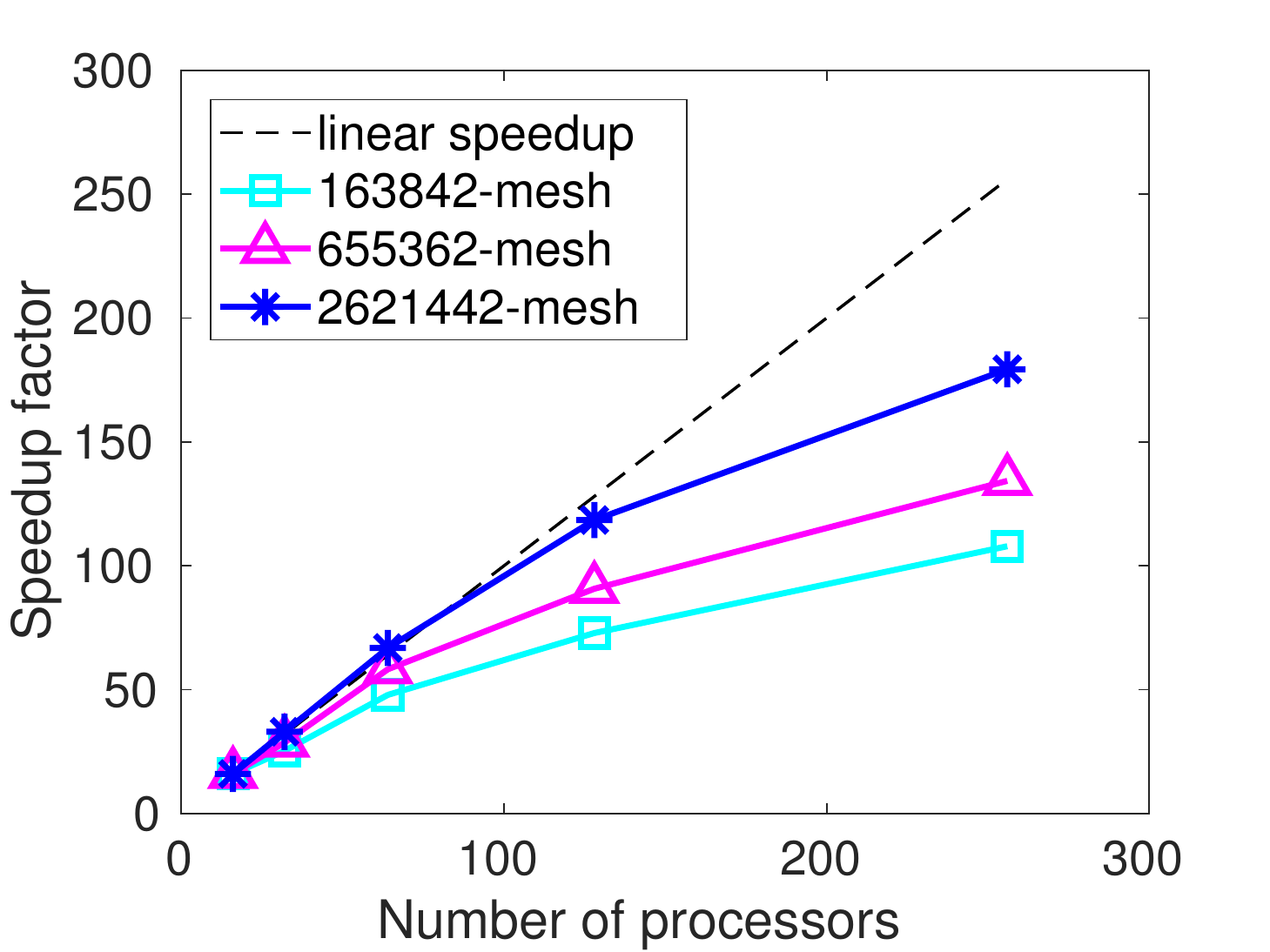}
\includegraphics[scale=0.36]{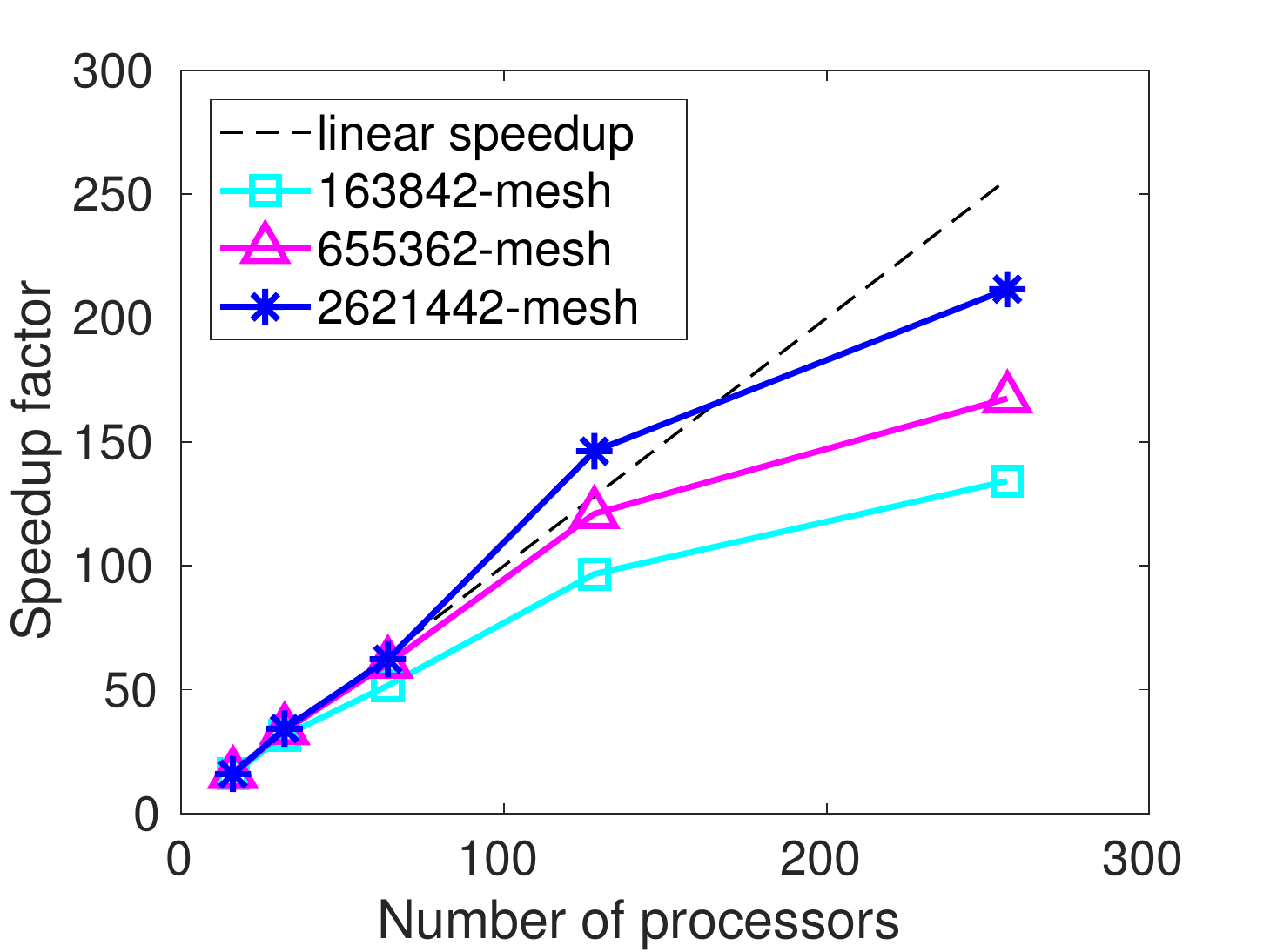}
\end{center}
\caption{Speedups of the parallel Lloyd P-LBFGS under  different number of  generators. From left to right: X3, X16, X64 densities.}
\label{scalingFactor}
\end{figure}

We first present in Figure~\ref{scalingFactor} the speedup of parallelization with varying number of generators (163,842, \;655,362, \;2,621,442) and different point-density 
functions. The speedup is measured using the cost per iteration.
It is observed that in generating relatively small meshes, for instance 163,842 generators, communication seems to dominate the overall execution time. That is why the 
scalability ends up being sub-linear. As the number of generators increases, the speedup gets better and better. When the number of grid points is in million scale, the 
parallel solver features almost linear scaling up to 256 processors.

\begin{figure}[h!]
\begin{center}
\includegraphics[scale=0.35]{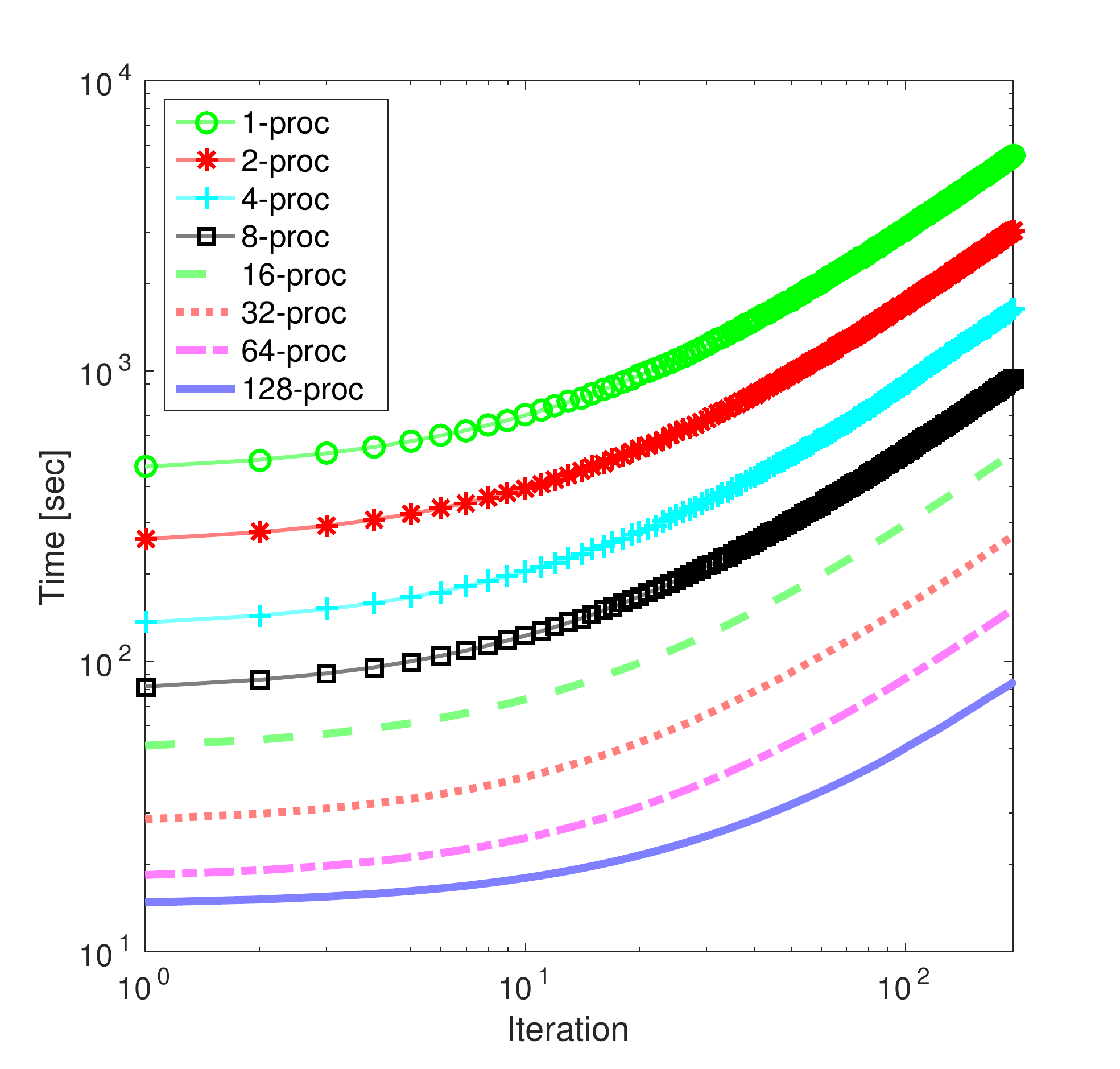}
\qquad
\includegraphics[scale=0.35]{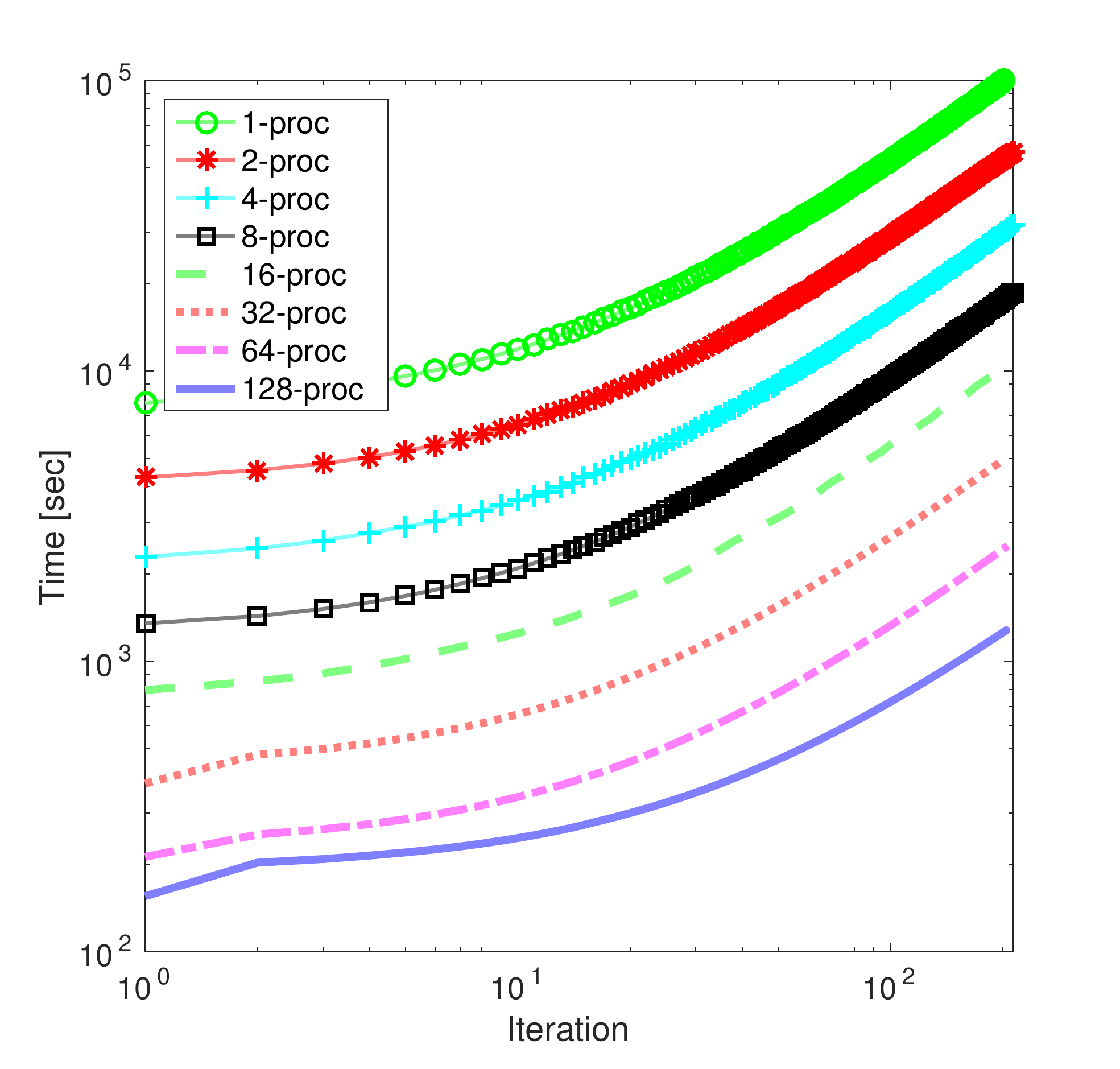}
\caption{Plots of running time vs. \# of iterations of the parallel Lloyd P-LBFGS with the X3 point-density. 
Left: 40,962 generators; right: 655,362 generators.}
\label{scalingTiter}
\end{center}
\end{figure}

\begin{figure}[h!]
\begin{center}
\includegraphics[scale=0.35]{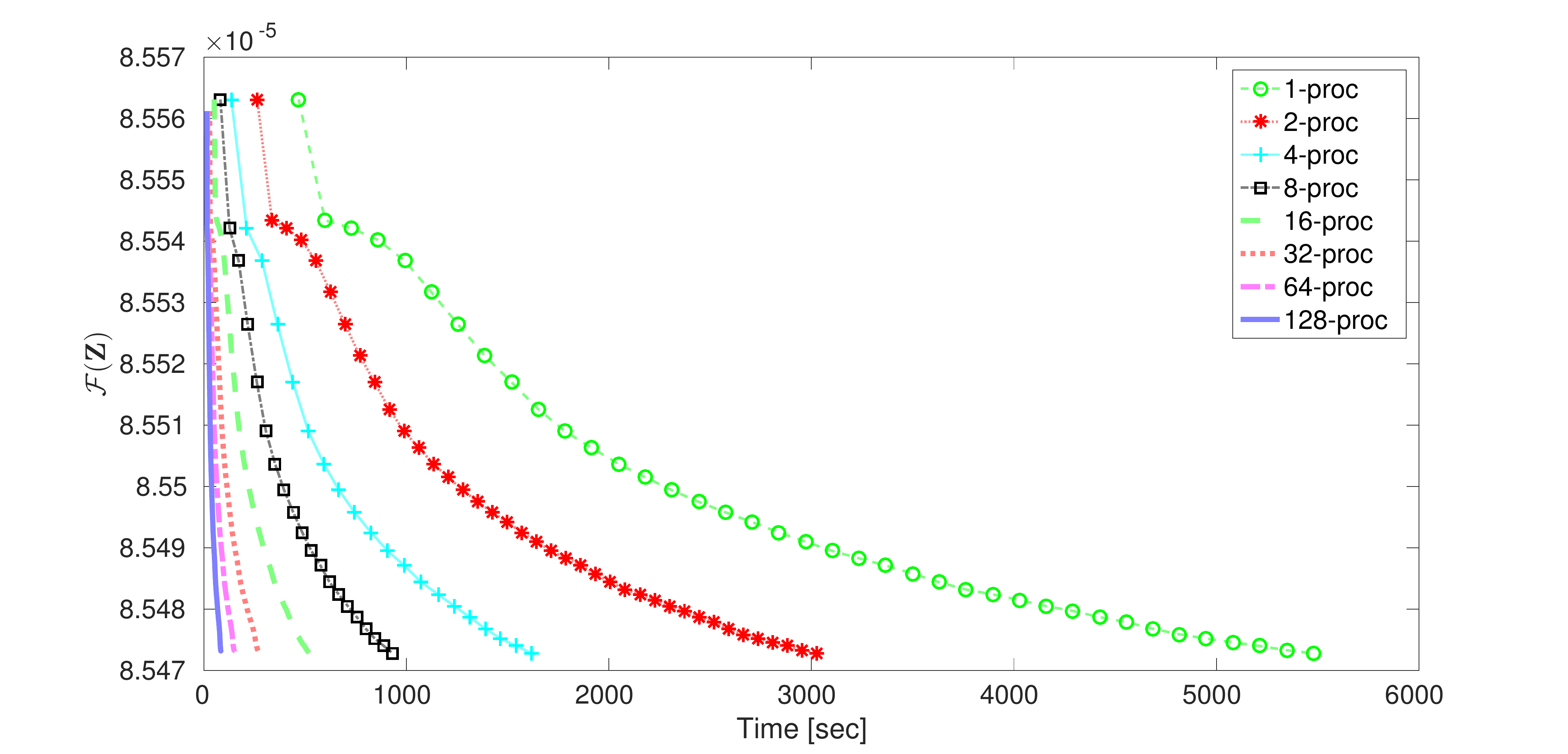}
\caption{Strong scaling efficiency check of the parallel Lloyd P-LBFGS  with 40,962 generators using X3 point-density: plot of $\mathcal{F}$-vs. running time.}
\label{scalingFT}
\end{center}
\end{figure}

We next show the scaling results towards generating the SCVTs with 655,362 generators. Initialized from a convergent CVT mesh of 2,562 generators, our code consecutively bisects and optimizes the mesh until a converged SCVT with 655,362 generators is obtained. The iteration history is recorded for the 655,362 generator mesh and also the 40,962 generator mesh generated along the consecutive bisection procedure. To measure the parallel efficiency, we launch the program with different number of processors ranging from 1 to 128. Figure~\ref{scalingTiter} shows the log-plot of the running time $t$ against the iteration number $n$. As they are related by $t=a\cdot n+c$, when $n$ is relatively large we can approximate $\log t \approx \log a + \log n$. Thus, the distance between the right ends of adjacent curves reflects the gain on time efficiency per iteration when the number of processors is doubled. Excellent scalability is confirmed because those distances are fairly evenly distributed. The plots of $\mathcal{F}$ and $\lVert\nabla\mathcal{F}\rVert$ vs. the running time for the 40,962 points mesh are shown in Figure~\ref{scalingFT}.

\begin{figure}[h!]
\begin{center}
\includegraphics[scale=0.5]{meshQ_bar2nd.png}\\
\includegraphics[scale=0.35]{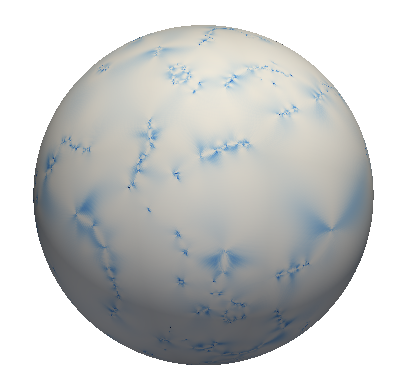}
\includegraphics[scale=0.35]{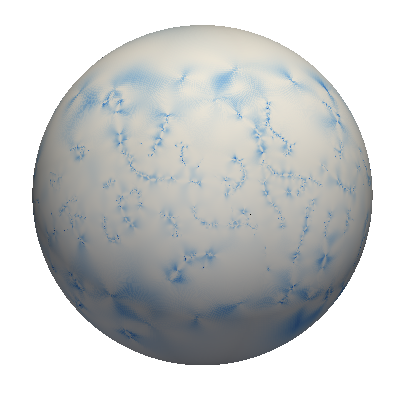}
\includegraphics[scale=0.35]{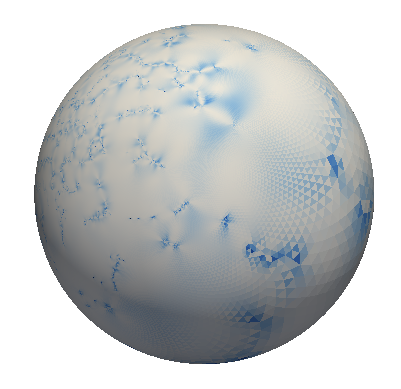}
\caption{Distribution  of triangle qualities of the  SCVT meshes with 655,362 generators created by the parallel Lloyd P-LBFGS using 128 processors. 
From left to right: X3, X16, X64 densities.}
\label{TriQpara}
\end{center}
\end{figure}

Finally, in Figure~\ref{TriQpara}, we plot the distribution of triangle qualities  of  the SCVT meshes with 655,362 generators  created by the parallel Lloyd P-LBFGS method with 128 processors. The total running times are  respectively 1296.4 seconds for X3 point-density, 1939.6 seconds for X16 point-density, and 1914.9 seconds for X64 point-density.

%%%%%%%%%%%%%%%%%%%%%%%%%%%%%%%%%%%%%%%%%%%%%%%%%%%%%%%%%%%%%%%%%%%%%%%%%%%
%%%%%%%%%%%%%%%%%%%%%%%%%%%%%%%%%%%%%%%%%%%%%%%%%%%%%%%%%%%%%%%%%%%%%%%%%%%
%%%%%%%%%%%%%%%%%%%%%%%%%%%%%%%%%%%%%%%%%%%%%%%%%%%%%%%%%%%%%%%%%%%%%%%%%%%
\section{Summary}\label{conclusion}

A Lloyd-preconditioned LBFGS method and {its} parallel implementation are proposed and numerically investigated for SCVT computation. Whereas LBFGS often speeds up the CVT computation as compared to the classical Lloyd method, it may lose ``optimal efficiency'' in generating  highly variable multi-resolution meshes. The Lloyd-preconditioned LBFGS method can further speeds up the CVT computation and  overcome this difficulty in large scale SCVT computation.

For quasi-uniform meshes, the Lloyd-preconditioned LBFGS method performs similarly as LBFGS, which is also the case when using the graph Laplacian preconditioner.  For typical multi-resolution meshes, however, the Lloyd-preconditioned LBFGS method significantly improves the performance of LBFGS, while also producing better mesh quality. Because the method enables efficient SCVT generation in highly variable multi-resolution, it is also useful for limited-area climate modeling by culling coarse grids. Note that it is often challenging to maintain regularity of meshes at boundaries by direct local mesh generation.

The parallel Lloyd-preconditioned LBFGS algorithm features well-balanced loading of grid points and performs well in the parallel {scalability} tests. It enables more convenient (faster mesh generation) and more
stable (better mesh quality) global climate modeling, especially for meshes with ultra high-resolution. 
The proposed  efficient CVT grid generator will thus play a fundamental role in next generations of climate modeling, where the spacing scale could be as small as 100 meters. Our interests will also include more accurate numerical simulations in ocean and atmosphere modeling based on SCVTs.

%Although the Lloyd-preconditioned LBFGS method is tested on sphere for climate modeling, the method is also applicable for general domains. A follow-up of this %
%work will be a more detailed application using data-driven point-density functions. Our interests will include more accurate simulations in ocean and atmosphere modeling.

%\section*{Acknowledgements}

\section*{References}

\bibliographystyle{model1b-num-names}

\begin{thebibliography}{}
\expandafter\ifx\csname url\endcsname\relax
  \def\url#1{\texttt{#1}}\fi
\expandafter\ifx\csname urlprefix\endcsname\relax\def\urlprefix{URL }\fi
\expandafter\ifx\csname href\endcsname\relax
  \def\href#1#2{#2} \def\path#1{#1}\fi

\end{thebibliography}


\begin{thebibliography}{24}
\expandafter\ifx\csname natexlab\endcsname\relax\def\natexlab#1{#1}\fi
\providecommand{\url}[1]{\texttt{#1}}
\providecommand{\href}[2]{#2}
\providecommand{\path}[1]{#1}
\providecommand{\DOIprefix}{doi:}
\providecommand{\ArXivprefix}{arXiv:}
\providecommand{\URLprefix}{URL: }
\providecommand{\Pubmedprefix}{pmid:}
\providecommand{\doi}[1]{\href{http://dx.doi.org/#1}{\path{#1}}}
\providecommand{\Pubmed}[1]{\href{pmid:#1}{\path{#1}}}
\providecommand{\bibinfo}[2]{#2}
\ifx\xfnm\relax \def\xfnm[#1]{\unskip,\space#1}\fi

%Type = Article
\bibitem[{Benzi and Tuma(1999)}]{benzi1999comparative}
\bibinfo{author}{M.~Benzi}, \bibinfo{author}{M.~Tuma}, \bibinfo{title}{A
  comparative study of sparse approximate inverse preconditioners},
  \bibinfo{journal}{Applied Numerical Mathematics} \bibinfo{volume}{30}
  (\bibinfo{year}{1999}) \bibinfo{pages}{305--340}.
  
%Type = Article
\bibitem[{Cignoni et~al.(1998)Cignoni, Montani and
  Scopigno}]{cignoni1998dewall}
\bibinfo{author}{P.~Cignoni}, \bibinfo{author}{C.~Montani},
  \bibinfo{author}{R.~Scopigno}, \bibinfo{title}{Dewall: A fast divide and
  conquer {D}elaunay triangulation algorithm in ed},
  \bibinfo{journal}{Computer-Aided Design} \bibinfo{volume}{30}
  (\bibinfo{year}{1998}) \bibinfo{pages}{333--341}.
  
%Type = Article
\bibitem[{Du and Emelianenko(2006)}]{DuAccel2006}
\bibinfo{author}{Q.~Du}, \bibinfo{author}{M.~Emelianenko},
  \bibinfo{title}{Acceleration schemes for computing centroidal {V}oronoi
  tessellations}, \bibinfo{journal}{Numer. Linear Algebra Appl.}
  \bibinfo{volume}{13} (\bibinfo{year}{2006}) \bibinfo{pages}{173--192}.
  
%Type = Article
\bibitem[{Du et~al.(2006)Du, Emelianenko and Ju}]{du2006convergence}
\bibinfo{author}{Q.~Du}, \bibinfo{author}{M.~Emelianenko},
  \bibinfo{author}{L.~Ju}, \bibinfo{title}{Convergence of the lloyd algorithm
  for computing centroidal {V}oronoi tessellations}, \bibinfo{journal}{SIAM
  journal on numerical analysis} \bibinfo{volume}{44} (\bibinfo{year}{2006})
  \bibinfo{pages}{102--119}.
  
%Type = Article
\bibitem[{Du et~al.(1999)Du, Faber and Gunzburger}]{du1999centroidal}
\bibinfo{author}{Q.~Du}, \bibinfo{author}{V.~Faber},
  \bibinfo{author}{M.~Gunzburger}, \bibinfo{title}{Centroidal {V}oronoi
  tessellations: Applications and algorithms}, \bibinfo{journal}{SIAM review}
  \bibinfo{volume}{41} (\bibinfo{year}{1999}) \bibinfo{pages}{637--676}.
  
%Type = Article
\bibitem[{Du and Gunzburger(2002)}]{du2002grid}
\bibinfo{author}{Q.~Du}, \bibinfo{author}{M.~Gunzburger}, \bibinfo{title}{Grid
  generation and optimization based on centroidal {V}oronoi tessellations},
  \bibinfo{journal}{Applied Mathematics and Computation} \bibinfo{volume}{133}
  (\bibinfo{year}{2002}) \bibinfo{pages}{591--607}.
  
%Type = Article
\bibitem[{Du et~al.(2003)Du, Gunzburger and Ju}]{du2003constrained}
\bibinfo{author}{Q.~Du}, \bibinfo{author}{M.~Gunzburger},
  \bibinfo{author}{L.~Ju}, \bibinfo{title}{Constrained centroidal {V}oronoi
  tessellations for surfaces}, \bibinfo{journal}{SIAM Journal on Scientific
  Computing} \bibinfo{volume}{24} (\bibinfo{year}{2003})
  \bibinfo{pages}{1488--1506}.
  
%Type = Article
\bibitem[{Dymond et~al.(2001)Dymond, Zhou and Deng}]{dymond2001}
\bibinfo{author}{P.~Dymond}, \bibinfo{author}{J.~Zhou},
  \bibinfo{author}{X.~Deng}, \bibinfo{title}{A 2-D parallel convex hull
  algorithm with optimal communication phases}, \bibinfo{journal}{Parallel
  Computing} \bibinfo{volume}{27} (\bibinfo{year}{2001})
  \bibinfo{pages}{243--255}.
  
%Type = Article
\bibitem[{Hateley et~al.(2015)Hateley, Wei and Chen}]{Hateley2015}
\bibinfo{author}{J.~Hateley}, \bibinfo{author}{H.~Wei},
  \bibinfo{author}{L.~Chen}, \bibinfo{title}{Fast methods for computing
  centroidal {V}oronoi tessellations}, \bibinfo{journal}{Journal of Scientific
  Computing} \bibinfo{volume}{63} (\bibinfo{year}{2015})
  \bibinfo{pages}{185--212}. 
%Type = Incollection

\bibitem[{Iri et~al.(1984)Iri, Murota and Ohya}]{iri1984fast}
\bibinfo{author}{M.~Iri}, \bibinfo{author}{K.~Murota},
  \bibinfo{author}{T.~Ohya}, \bibinfo{title}{A fast {V}oronoi-diagram algorithm
  with applications to geographical optimization problems}, in:
  \bibinfo{booktitle}{System Modelling and Optimization},
  \bibinfo{publisher}{Springer}, \bibinfo{year}{1984}, pp.
  \bibinfo{pages}{273--288}.
  
%Type = Phdthesis
\bibitem[{Jacobsen(2011)}]{JacobThesis}
\bibinfo{author}{D.~Jacobsen}, \bibinfo{title}{Parallel Grid Generation and
  Multi-Resolution Methods for Climate Modeling Applications}, Ph.D. thesis,
  Florida State University, \bibinfo{year}{2011}.
  
%Type = Article
\bibitem[{Jacobsen et~al.(2013)Jacobsen, Gunzburger, Ringler, Burkardt and
  Peterson}]{JacobMPI2013}
\bibinfo{author}{D.~Jacobsen}, \bibinfo{author}{M.~Gunzburger},
  \bibinfo{author}{T.~Ringler}, \bibinfo{author}{J.~Burkardt},
  \bibinfo{author}{J.~Peterson}, \bibinfo{title}{Parallel algorithms for planar
  and spherical {D}elaunay construction with an application to centroidal
  {V}oronoi tessellations}, \bibinfo{journal}{Geoscientific Model Development}
  \bibinfo{volume}{6} (\bibinfo{year}{2013}) \bibinfo{pages}{1353--1365}.
  
%Type = Techreport
\bibitem[{Jiang et~al.(2004)Jiang, Byrd, Eskow and
  Schnabel}]{jiang2004preconditioned}
\bibinfo{author}{L.~Jiang}, \bibinfo{author}{R.~Byrd},
  \bibinfo{author}{E.~Eskow}, \bibinfo{author}{R.~Schnabel}, \bibinfo{title}{A
  preconditioned {L-BFGS} algorithm with application to molecular energy
  minimization}, \bibinfo{type}{Technical Report}, DTIC Document,
  \bibinfo{year}{2004}.
  
%Type = Article
\bibitem[{Ju et~al.(2002)Ju, Du and Gunzburger}]{ju2002probabilistic}
\bibinfo{author}{L.~Ju}, \bibinfo{author}{Q.~Du},
  \bibinfo{author}{M.~Gunzburger}, \bibinfo{title}{Probabilistic methods for
  centroidal {V}oronoi tessellations and their parallel implementations},
  \bibinfo{journal}{Parallel Computing} \bibinfo{volume}{28}
  (\bibinfo{year}{2002}) \bibinfo{pages}{1477--1500}.

%Type = Article
\bibitem[{Liu et~al.(2009)Liu, Wang, L{\'e}vy, Sun, Yan, Lu and
  Yang}]{Liu2009CVT}
\bibinfo{author}{Y.~Liu}, \bibinfo{author}{W.~Wang},
  \bibinfo{author}{B.~L{\'e}vy}, \bibinfo{author}{F.~Sun},
  \bibinfo{author}{D.~Yan}, \bibinfo{author}{L.~Lu}, \bibinfo{author}{C.~Yang},
  \bibinfo{title}{On centroidal {V}oronoi tessellation: Energy smoothness and
  fast computation}, \bibinfo{journal}{ACM Transactions on Graphics}
  \bibinfo{volume}{28} (\bibinfo{year}{2009}) \bibinfo{pages}{101}.
%Type = Article
\bibitem[{Lloyd(1982)}]{lloyd1982least}
\bibinfo{author}{S.~Lloyd}, \bibinfo{title}{Least squares quantization in pcm},
  \bibinfo{journal}{IEEE transactions on information theory}
  \bibinfo{volume}{28} (\bibinfo{year}{1982}) \bibinfo{pages}{129--137}.
  
%Type = Inproceedings
\bibitem[{MacQueen et~al.(1967)}]{macqueen1967some}
\bibinfo{author}{J.~MacQueen}, et~al., \bibinfo{title}{Some methods for
  classification and analysis of multivariate observations}, in:
  \bibinfo{booktitle}{Proceedings of the fifth Berkeley symposium on
  mathematical statistics and probability}, volume~\bibinfo{volume}{1},
  \bibinfo{organization}{Oakland, CA, USA}, pp. \bibinfo{pages}{281--297}.
  
%Type = Book
\bibitem[{Nocedal and Wright(2006)}]{Nocedal2006NO}
\bibinfo{author}{J.~Nocedal}, \bibinfo{author}{S.J. Wright},
  \bibinfo{title}{Numerical Optimization}, \bibinfo{edition}{2nd} ed.,
  \bibinfo{publisher}{Springer}, \bibinfo{address}{New York},
  \bibinfo{year}{2006}.
  
%Type = Article
\bibitem[{Ringler et~al.(2011)Ringler, Jacobsen, Gunzburger, Ju, Duda and
  Skamarock}]{ringler2011exploring}
\bibinfo{author}{T.~Ringler}, \bibinfo{author}{D.~Jacobsen},
  \bibinfo{author}{M.~Gunzburger}, \bibinfo{author}{L.~Ju},
  \bibinfo{author}{M.~Duda}, \bibinfo{author}{W.~Skamarock},
  \bibinfo{title}{Exploring a multiresolution modeling approach within the
  shallow-water equations}, \bibinfo{journal}{Monthly Weather Review}
  \bibinfo{volume}{139} (\bibinfo{year}{2011}) \bibinfo{pages}{3348--3368}.
  
%Type = Article
\bibitem[{Ringler et~al.(2013)Ringler, Petersen, Higdon, Jacobsen, Jones and
  Maltrud}]{ringler2013multi}
\bibinfo{author}{T.~Ringler}, \bibinfo{author}{M.~Petersen},
  \bibinfo{author}{R.~Higdon}, \bibinfo{author}{D.~Jacobsen},
  \bibinfo{author}{P.~Jones}, \bibinfo{author}{M.~Maltrud}, \bibinfo{title}{A
  multi-resolution approach to global ocean modeling}, \bibinfo{journal}{Ocean
  Modelling} \bibinfo{volume}{69} (\bibinfo{year}{2013})
  \bibinfo{pages}{211--232}.
  
%Type = Article
\bibitem[{Sakaguchi et~al.(2015)Sakaguchi, Leung, Zhao, Yang, Lu, Hagos,
  Rauscher, Dong, Ringler and Lauritzen}]{Sakaguchi2015}
\bibinfo{author}{K.~Sakaguchi}, \bibinfo{author}{L.~Leung},
  \bibinfo{author}{C.~Zhao}, \bibinfo{author}{Q.~Yang},
  \bibinfo{author}{J.~Lu}, \bibinfo{author}{S.~Hagos},
  \bibinfo{author}{S.~Rauscher}, \bibinfo{author}{L.~Dong},
  \bibinfo{author}{T.~Ringler}, \bibinfo{author}{P.~Lauritzen},
  \bibinfo{title}{Exploring a multiresolution approach using {AMIP}
  simulations}, \bibinfo{journal}{Journal of Climate} \bibinfo{volume}{28}
  (\bibinfo{year}{2015}) \bibinfo{pages}{5549--5574}.

%Type = Article
\bibitem[{Sakaguchi et~al.(2016)Sakaguchi, Lu, Leung, Zhao, Li and
  Hagos}]{JAME20337}
\bibinfo{author}{K.~Sakaguchi}, \bibinfo{author}{J.~Lu},
  \bibinfo{author}{L.~Leung}, \bibinfo{author}{C.~Zhao},
  \bibinfo{author}{Y.~Li}, \bibinfo{author}{S.~Hagos}, \bibinfo{title}{Sources
  and pathways of the upscale effects on the southern hemisphere jet in
  {MPAS-CAM4} variable-resolution simulations}, \bibinfo{journal}{Journal of
  Advances in Modeling Earth Systems} \bibinfo{volume}{8}
  (\bibinfo{year}{2016}) \bibinfo{pages}{1786--1805}.

%Type = Incollection
\bibitem[{Shewchuk(1996)}]{Shewchuk1996}
\bibinfo{author}{J.~Shewchuk}, \bibinfo{title}{Triangle: Engineering a {2D}
  quality mesh generator and {D}elaunay triangulator}, in:
  \bibinfo{booktitle}{Applied Computational Geometry Towards Geometric
  Engineering}, \bibinfo{publisher}{Springer Berlin Heidelberg},
  \bibinfo{year}{1996}, pp. \bibinfo{pages}{203--222}.
  
%Type = Article
\bibitem[{Zhao et~al.(2016)Zhao, Leung, Park, Hagos, Lu, Sakaguchi, Yoon,
  Harrop, Skamarock and Duda}]{JAME20336}
\bibinfo{author}{C.~Zhao}, \bibinfo{author}{L.~Leung},
  \bibinfo{author}{S.~Park}, \bibinfo{author}{S.~Hagos},
  \bibinfo{author}{J.~Lu}, \bibinfo{author}{K.~Sakaguchi},
  \bibinfo{author}{J.~Yoon}, \bibinfo{author}{B.~Harrop},
  \bibinfo{author}{W.~Skamarock}, \bibinfo{author}{M.~Duda},
  \bibinfo{title}{Exploring the impacts of physics and resolution on
  aqua-planet simulations from a nonhydrostatic global variable-resolution
  modeling framework}, \bibinfo{journal}{Journal of Advances in Modeling Earth
  Systems} \bibinfo{volume}{8} (\bibinfo{year}{2016})
  \bibinfo{pages}{1751--1768}. 
\end{thebibliography}

\end{document}